\documentclass[11pt]{article}

\usepackage[small,bf]{caption}
\usepackage{amsmath, amsthm, amssymb}
\usepackage{mathtools}
\usepackage{algorithmic}
\usepackage{algorithm}
\usepackage[parfill]{parskip}
\usepackage{soul}
\usepackage{multirow}
\usepackage{caption}
\usepackage{subcaption}

\usepackage{epsfig,fullpage,psfrag,url,verbatim}
\usepackage{graphics}
\usepackage{graphicx}
\usepackage{epstopdf}

\usepackage{amsfonts}
\usepackage{color}

\usepackage{array}
\usepackage{booktabs}
\usepackage[dvipsnames]{xcolor}
\usepackage{tikz}

\newcolumntype{M}[1]{>{\centering\arraybackslash}m{#1}}

\newtheorem{cor}{Corollary}[section]
\newtheorem{thm}{Theorem}[section]
\newtheorem{lem}{Lemma}[section]
\newtheorem{prop}{Proposition}[section]

\newtheorem{mydef}{Definition}[section]
\newtheorem{rem}{Remark}[section]
\newtheorem{ass}{Assumption}[section]

\newtheorem{exam}{Example}[section]

\newcommand{\dl}{\dot{l}}

\newcommand{\dom}{\mathrm{dom}}

\newcommand{\hh}{h(\cdot)}
\newcommand{\ff}{f(\cdot)}
\newcommand{\EE}{\mathbb{E}}

\newcommand{\Dh}{D_h}

\newcommand{\tg}{\tilde{g}}

\newcommand{\RR}{\mathbb{R}}
\newcommand{\bbeta}{\bar{\beta}}
\newcommand{\tbeta}{\tilde{\beta}}
\newcommand{\dlsji}{\dot l^*_{j_i}}
\newcommand{\dls}{\dot l^*}
\newcommand{\wiji}{w^i_{j_i}}

\newcommand{\hw}{\hat w}
\newcommand{\hbeta}{\hat \beta}
\newcommand{\bw}{\bar w}
\newcommand{\ls}{l^*}
\newcommand{\bI}{{\bf I}}

\DeclareMathOperator*{\argmax}{arg\,max}

\begin{document}

\title{Generalized Stochastic Frank-Wolfe Algorithm with Stochastic ``Substitute'' Gradient for Structured Convex Optimization}

\author{Haihao Lu\thanks{MIT Department of Mathematics, 77 Massachusetts Avenue, Cambridge, MA   02139
({mailto:  haihao@mit.edu}).}
\and
Robert M. Freund\thanks{MIT Sloan School of Management, 77 Massachusetts Avenue, Cambridge, MA   02139
({mailto:  rfreund@mit.edu}).  This author's research is supported by AFOSR Grant No. FA9550-15-1-0276.}
}


\maketitle

\begin{abstract}
The stochastic Frank-Wolfe method has recently attracted much general interest in the context of optimization for statistical and machine learning due to its ability to work with a more general feasible region. However, there has been a complexity gap \textcolor{black}{in the dependence on the optimality tolerance $\varepsilon$} in the guaranteed convergence rate for stochastic Frank-Wolfe compared to its deterministic counterpart.  In this work, we present a new generalized stochastic Frank-Wolfe method which closes this gap for the class of structured optimization problems encountered in statistical and machine learning characterized by empirical loss minimization with a certain type of ``linear prediction'' property (formally defined in the paper), which is typically present in loss minimization problems in practice.  Our method also introduces the notion of a ``substitute gradient" that is a not-necessarily-unbiased sample of the gradient.  We show that our new method is equivalent to a particular randomized coordinate mirror descent algorithm applied to the dual problem, which in turn provides a new interpretation of randomized dual coordinate descent in the primal space.  Also, in the special case of a strongly convex regularizer our generalized stochastic Frank-Wolfe method (as well as the randomized dual coordinate descent method) exhibits linear convergence.  Furthermore, we present computational experiments that indicate that our method outperforms other stochastic Frank-Wolfe methods \textcolor{black}{for a sufficiently small optimality tolerance}, consistent with the theory developed herein.

\end{abstract}

\section{Introduction}\label{intro}

Our problem of interest is the following optimization problem:

\begin{equation}\label{eq:primal}
\mathrm{(P):} \ \ \ \ \ \ \ \min_\beta P(\beta):=\tfrac{1}{n }\sum_{j=1}^{n}l_{j}(x_j^T\beta)+R(\beta) \ ,
\end{equation}
where $\beta \in \RR^p$, $l_j(\cdot) : \RR \rightarrow \RR$, $j=1, \ldots, n$, is a univariate function (the $j^{\mathrm{th}}$ loss function), $s_j = x_j^T\beta$ is the ``predicted value'' of the model $\beta$ for the data sample $x_j$, and $R(\cdot)$ is some other function that can be used to model a variety of properties such as a regularizer, an indicator function of a feasible region $Q$, a penalty term, coupling constraints, etc.  Notice that the predicted value $s_j= x_j^T\beta$ is a linear function of $\beta$.  Hence in the intentional context of empirical loss minimization we refer to \eqref{eq:primal} as (empirical) loss minimization with ``linear prediction.''  This structure arises quite generally throughout statistical and machine learning, which we illustrate below via several salient examples.  Throughout this paper we assume the following regarding the functions in our problem setup \eqref{eq:primal}:\medskip

\begin{ass}\label{ass:smooth}  The following hold:
\begin{enumerate}
\item for $j=1, \ldots, n$, the univariate function $l_j(\cdot)$ is $\gamma$-smooth, namely $| \dot l_j(a) - \dot l_j(b)| \le \gamma |a-b|$ for all $a$, $b$, and is strictly convex,
\item $\dom R(\cdot)$ is bounded, and the subproblem
\begin{equation}\label{eq:sub}
\min_\beta c^T\beta +R(\beta)
\end{equation}
attains its optimum and can be easily solved for any $c$, and
\item $ 0 \in \dom R(\cdot)$.
\end{enumerate}
\end{ass}

We note regarding {\it (1.)} above that strict convexity (instead of simple convexity) is only needed to guarantee that the conjugate function $l_j^*(\cdot)$ is differentiable.  Regarding {\it (2.)}, this is a generalization of a linear optimization oracle as follows:  in the case when $R(\cdot)$ is the indicator function $\bI_Q(\cdot)$ of a set $Q \subset \RR^p$ (namely, $\bI_Q(\beta):=0$ if $\beta \in Q$, and $\bI_Q(\beta) :=+\infty$ otherwise), then $Q$ is the feasible region of (P), and {\it (2.)} states that the feasible region $Q$ is bounded and that it is easy to solve linear optimization problems on $Q$.  Also, {\it(3.)} above is for notational convenience, as we can always translate a given feasible point so that $0 \in \dom R(\cdot)$.

\subsection{Examples in Statistical and Machine Learning}
Here we present several applications of our problem setup \eqref{eq:primal} in statistical and machine learning.  (For other applications particularly amenable to solution by the Frank-Wolfe method, we refer the reader to \cite{jaggi2013revisiting}.)\medskip

\begin{exam} {\bf LASSO\cite{tibshirani1996regression}, ridge regression\cite{hoerl1970ridge}, sparse logisitic regression\cite{ravikumar2010high}.} Consider the least-squares regression problem where a set of training samples $\{(x_j, y_j)\}_{j=1}^n$ is given. The LASSO optimization problem (in constraint format) is:
\begin{equation*}
\begin{array}{cl}
\min_{\beta} & \tfrac{1}{2n}\sum_{j=1}^n (y_j-x_j^T \beta)^2\\ \\
s.t. & \|\beta\|_1 \le \delta \ ,
\end{array}
\end{equation*}
which is an instance of (P) by using the least squares loss function $l_j(\cdot)=\frac{1}{2}(y_j - \cdot)^2$ and using the indicator function of an $\ell_1$ ball as the ``regularizer'' function $R(\cdot)$, namely $R(\beta):=\bI_{\{\|\beta\|_1\le \delta \}}(\beta)$.

The ridge regression optimization problem adds the regularizer $\tfrac{\lambda}{2} \|\beta\|_2^2$ to the least squares objective function for the parameter $\lambda >0$, and omits the $\ell_1$ ball constraint.  Notice that because $\beta = 0$ is a feasible solution it follows that the optimal objective value is bounded above by $ \|y\|^2_2/(2n)$, and therefore we can model $R(\beta)=\frac{\lambda}{2}\|\beta\|_2^2 + \bI_{\{\|\beta\|^2_2 \le \|y\|^2_2 /(n\lambda) \}}(\beta)$, which ensures that $\dom R(\cdot)$ is bounded.

The $\ell_1$-regularized logistic regression optimization problem seeks a solution of:
\begin{equation*}
\min_{\beta} P(\beta) = \tfrac{1}{n}\sum_{j=1}^n \ln (1 + \exp(-y_jx_j^T \beta)) + \lambda\|\beta\|_1\ ,
\end{equation*}
for a given set of training samples $\{(x_j, y_j)\}_{j=1}^n$ where $y_j \in \{-1,1\}$, and is an instance of (P) using the logistic loss functions  $l_j(\cdot)=\ln(1+\exp(-y_j\cdot))$ with $R(\beta)=\lambda\|\beta\|_1 + \bI_{\left\{\|\beta\|_1 \le \ln(2)/\lambda \right\}}(\beta)$ where the indicator function term is structurally redundant but is added as in the previous example to ensure that $\dom R(\cdot)$ is bounded.
\end{exam}\medskip

\begin{exam} {\bf Matrix completion\cite{fazel2002matrix}\cite{candes2009exact}.}\label{robespierre} In the matrix completion problem, we seek to compute a low-rank matrix that well-approximates a given matrix $M \in \RR^{n \times p}$ on the set $\Omega$ of observed entries $(i,j)$.  The convex relaxation of this problem is the following nuclear-norm optimization problem:
\begin{equation*}
\begin{array}{ll}
\min_{\beta \in \RR^{n \times p}} & \frac{1}{2|\Omega|}\sum_{(i,j)\in \Omega} (M_{i,j} - \beta_{i,j})^2 \\ \\
s.t. & \|\beta\|_* \le \delta \ ,
\end{array}
\end{equation*}
where $\|\cdot\|_*$ is the nuclear norm. In order to translate the matrix completion problem to the setting of (P), we consider any index pair $(i,j)\in \Omega$ as a sample, and we have $l_{(i,j)}(\cdot)=\tfrac{1}{2}(\cdot-M_{i,j})^2$,  and $R(\beta)= \bI_{\{ \|\beta\|_* \le \delta \}}(\beta)$.
\end{exam}\medskip

\begin{exam} {\bf Structured sparse matrix estimation with CUR factorization\cite{mahoney2009cur}\cite{mairal2011convex}.} We seek to compute an approximate factorization $M\approx CUR$ of a given data matrix $M\in \RR^{n\times d}$ such that $C$ contains a subset of $c$ columns from $M$ and $R$ contains a subset of $r$ rows from $M$. Mairal et al. \cite{mairal2011convex} proposed the following convex relaxation of this problem:
\begin{equation*}
\begin{array}{cl}
\min_{\beta} & \tfrac{1}{2nd}\|M-M\beta M\|_F^2 \\ \\
s.t. & \sum_{i}\|\beta_{i, \cdot}\|_\infty \le \delta \\ \\
& \sum_{j}\|\beta_{\cdot, j}\|_\infty \le \delta \ ,
\end{array}
\end{equation*}
which is an instance of (P) by modeling the $(i,j)^{\mathrm{th}}$ loss term in (P) as $\frac{1}{2}(M_{i,j}-M_i^T \beta M_j)^2$ (which is a least squares loss of a particular linear function of the matrix variable $\beta$), and $R(\beta)=\bI_{\{\sum_{i}\|\beta_{i,\cdot}\|_\infty \le \delta, \ \sum_{j}\|\beta_{\cdot, j}\|_\infty \le \delta\}}(\beta)$.
\end{exam}

\subsection{Stochastic Generalized Frank-Wolfe with Stochastic Substitute Gradient}

Let $X \in \RR^{n \times p}$ denote the data matrix whose rows are comprised of the vectors $x_1, \ldots, x_n$, i.e., the $j^{\mathrm{th}}$ row of $X$ is the vector $x_j$, $j=1, \ldots, n$.  Let us define $L(s): \RR^n \rightarrow \RR$ by $L(s) :=\sum_{j=1}^n l_j(s_j)$ which is the total losses associated with $s \in \RR^n$.  We refer to $s=X\beta$ as the ``predicted values'' in the context of empirical loss minimization.

Algorithm \ref{al:SFW} presents the main algorithmic contribution of this paper, which is a first-order method for tackling the problem (P).   We call the method ``Stochastic Generalized Frank-Wolfe method with Stochastic Substitute Gradient" for reasons which we will discuss as we walk through the structure of the method below.

\begin{algorithm}
\caption{Stochastic Generalized Frank-Wolfe with Stochastic Substitute Gradient}\label{al:SFW}
$ \ $
\begin{algorithmic}
\STATE {\bf Initialize.}  Initialize with $\bbeta^{-1}=0$, $s^0=0$, and substitute gradient $d^0=\frac{1}{n}X^T\nabla L(s^0)$, with
step-size sequences $\{\alpha_{i}\} \in (0,1]$ and $\{\eta_{i}\} \in (0,1]$.  \medskip

\STATE  For iterations $i=0,1, \ldots $

\ \ \ \ \ \ \  {\bf Solve l.o.o. subproblem:} Compute $\tbeta^{i} \in \arg\min_{\beta}\left\{\left( d^i\right)^{T}\beta+R(\beta)\right\}$

\ \ \ \ \ \ \  {\bf Choose random index:}  Choose $j_{i} \in {\cal U}[1, \ldots, n]$

\ \ \ \ \ \ \  {\bf Update $s$ value:} $s_{j_i}^{i+1}\leftarrow(1-\eta_{i})s_{j_i}^{i}+\eta_{i}(x_{j_{i}}^{T}\tbeta^{i})$,  and $s_{j}^{i+1} \gets s_{j}^{i}$ for $j \ne j_i$ \medskip

\ \ \ \ \ \ \  {\bf Update substitute gradient:}  $d^{i+1}= \frac{1}{n}X^T\nabla L(s^{i+1}) =  d^i+\frac{1}{n}\left(\dot{l}_{j_{i}}(s^{i+1}_{j_i})-\dot{l}_{j_{i}}(s^i_{j_{i}})\right)x_{j_{i}}$\\ \medskip

\ \ \ \ \ \ \  {\bf Update primal variable:} $\bbeta^{i} \leftarrow (1-\alpha_i) \bbeta^{i-1} + \alpha_i \tbeta^i$.\medskip

\ \ \ \ \ \ \  {\bf (Optional Accounting:)} $w^{i+1} \gets \nabla L(s^{i+1})$ \medskip
\end{algorithmic}
\end{algorithm}\medskip

We can write the first part of the objective function of (P) as $f(\beta):=\tfrac{1}{n} L(X\beta) = \tfrac{1}{n} L(s)$ with $s=X\beta$.  We have $\nabla L(s) =(\dl_1(s_1), \ldots, \dl_n(s_n))$ and the gradient of $\ff$ can be written as
\begin{equation}\label{rob1} \nabla f(\beta) = \tfrac{1}{n} X^T\nabla L(X\beta) = \tfrac{1}{n}\sum_{j=1}^n \dl_j(x_j^T\beta) x_j  \ , \end{equation}
which we can re-write as $ \nabla f(\beta) = \tfrac{1}{n} X^T w \ \mathrm{where} \ w = \nabla L(s) \ \mathrm{and} \ s = X\beta$, and which can be alternatively stated as:
\begin{equation}\label{rob2} \nabla f(\beta) = \tfrac{1}{n} \sum_{j=1}^n w_j x_j \ \mathrm{where} \ w_j = \dl_j(s_j) \ \mathrm{and} \ s_j = x_j^T\beta \ , \ j=1, \ldots, n \ . \end{equation} Here we emphasize that $w$ is the vector of weights on the data values $X$ in the composition of the gradient, and $s$ is the vector of predicted values $X\beta$.

Especially in the context of ``big data'' applications of statistical and machine learning where $n$ is huge, it can be extremely expensive to compute $\nabla \ff$.  We therefore maintain a ``substitute gradient'' in Algorithm \ref{al:SFW} that is constructed stochastically.  This is accomplished as follows:  let $\bbeta^{i-1}$ be the value of $\beta$ at the start of iteration $i$ of the method, and we have a substitute gradient $d^i$ that is the current proxy/substitute for the true gradient $\nabla f(\bbeta^{i-1})$, where $d^i$ is computed by:
\begin{equation}\label{rob3}
d^i = \tfrac{1}{n} \sum_{j=1}^n w^i_j x_j \ \mathrm{where} \ w^i_j = \dl_j(s^i_j)  \ , \ j=1, \ldots, n \ ,
\end{equation}
for a given $s^i$ that is the value of $s$ at iteration $i$.  But in contrast to \eqref{rob2} it will \underline{not} necessarily hold that $s^i_j = x_j^T\bbeta^i$ for $j=1, \ldots, n $ (equivalently $s^i = X\bbeta^i$).  (In fact, $d^i$ will not necessarily be an unbiased estimate of $\nabla f(\bbeta^i)$ as this will not be needed.)  In the identical spirit as randomized coordinate descent, $s^{i+1}$ will be determined by choosing a random index $j_i \in {\cal U}[1, \ldots, n]$ and updating only the coordinate $j_i$ of $s^i$, so that $s^{i+1} = s^i + \Delta^i e_{j_i}$ for some specific iteration-dependent scalar $\Delta^i$ (where $e_\ell$ denotes the $\ell^{\mathrm{th}}$ unit coordinate vector in $\RR^n$).  This is accomplished in the ``choose random index'' step and the ``update $s$ value'' step in Algorithm \ref{al:SFW}.

We now walk through the structure of Algorithm \ref{al:SFW} in complete detail.  The method is initialized with the initial decision variable $\beta$ set to $\bbeta^{-1} = 0$ and the vector of predicted values $s^0 = X\bbeta^{-1} = 0$ and initial substitute gradient $d^0=\frac{1}{n}X^T\nabla L(s^0)$, which corresponds to the true predicted value and true gradient at $\bbeta^{-1} = 0$. In iteration $i$, we use the substitute gradient $d^i$ to compute $\tbeta^i$, which is a solution to the (generalized) linear optimization oracle (``l.o.o."), where recall that this step specifies to solving a linear optimization problem over a set $Q$ in the specific case when the $R(\cdot)$ is the indicator function of $Q$, namely $R(\cdot) = \bI_{Q}(\cdot)$.  Regarding updating the current predicted values $s^i$, we randomly choose a sample (a coordinate) $j_i$ and only update $s_{j_i}$ as a certain convex combination of the current predicted value $s^i_{j_i}$ and the predicted value for the $j_i^{\mathrm{th}}$ sample at $\tbeta^i$, namely $x_{j_i}^T\tbeta^i$, so that $s^{i+1}_{j_i} \leftarrow(1-\eta_{i})s_{j_i}^{i}+\eta_{i}(x_{j_{i}}^{T}\tbeta^{i})$. Then we update the substitute gradient to make sure that $d^{i+1}=\frac{1}{n}X^T\nabla L(s^{i+1})$.   The last step at iteration $i$ is to take a Frank-Wolfe step to update $\bbeta^i \leftarrow (1-\alpha_i) \bbeta^{i-1} + \alpha_i \tbeta^i$ by taking a convex combination of the previous primal variable value $\bbeta^{i-1}$ and the solution $\tbeta^i$ of the just-solved linear optimization oracle.  Finally -- and ``optionally'' since it does not affect future computations -- we can perform an optional accounting step to update the dual variable $w^{i+1} \gets \nabla L(s^{i+1})$ in order to compute a duality gap certificate if desired.  (The nature of this duality will be understood once we look at the dual problem of (P) in Section \ref{dualdual}.)

Note that the computations in Algorithm \ref{al:SFW} are minimally affected by the dimension $n$.  Except for the initial computation of the gradient $d^0$ which is $O(np)$ operations, $s^i$ and $w^i$ are only updated by one coefficient at each iteration, and $d^{i+1}$ is updated by adding a scalar multiple of $x_{j_i}$ to $d^i$, which is $O(p)$ operations.  The updates of $\bbeta^i$ are $O(p)$ operations after solving for the optimal value $\tbeta^i$ in the linear optimization oracle, which is assumed to be easy to compute.

It is useful to place Algorithm \ref{al:SFW} in the context of the (deterministic) Frank-Wolfe method.  The Frank-Wolfe method is designed primarily to tackle the constrained convex optimization problem:  $\min_{\beta \in Q} f(\beta)$ where $\ff$ is a smooth convex function and $Q$ is a convex body, and it is assumed that linear optimization over $Q$ is easy to compute.  The optimization problem can of course be re-written as $\min_{\beta} f(\beta) + R(\beta)$ with $R(\cdot) = \bI_{Q}(\cdot)$.  The Frank-Wolfe update is:
\begin{equation}
\begin{array}{lcl}
\tbeta^i  \in   \arg\min_{\beta\in Q} \left\{\nabla f(\beta^i)^T \beta\right\}  \ \ \mathrm{and} \ \
\beta^{i+1} = (1-\alpha_i)\beta^i + \alpha_i \tbeta^i \ .
\end{array}
\end{equation}
It can be shown that with an appropriate choice of step-size sequence $\{\alpha_i\}$ that the Frank-Wolfe method computes an $\varepsilon$-optimal solution in $O\left(\tfrac{LD^2}{\varepsilon} \right)$ iterations, where $L$ is the Lipschitz constant of $f(\cdot)$ on $Q$ and $D$ is the diameter of $Q$, see \cite{frank-wolfe}, \cite{jaggi2013revisiting}, and \cite{freund2016new}.  Since our focus (for stochastic versions of Frank-Wolfe) will be on the dependence on $\varepsilon$, we will typically ignore these other instance-dependent constants and write the above as $O(\frac{1}{\varepsilon})$ iterations.

Due to its low iteration cost and convenient structural properties, the Frank-Wolfe method is especially applicable in several areas of statistical and machine learning and has thus received much renewed interest in recent years, see \cite{jaggi2013revisiting}, \cite{harchaoui2015conditional}, \cite{freund2017extended}, \cite{freund2016new}, and the references therein.  The Frank-Wolfe method can be generalized to deal with the more general problem $\min_\beta f(\beta) + R(\beta)$ where $R(\cdot)$ is any convex function with bounded domain and for which the ``linear optimization problem'' $\min_\beta c^T\beta + R(\beta)$ is easy to compute.  The generalized Frank-Wolfe update then is:
\begin{equation}
\begin{array}{lcl}
\tbeta^i  \in   \arg\min_{\beta} \left\{\nabla f(\beta^i)^T \beta + R(\beta) \right\}  \ \ \mathrm{and} \ \
\beta^{i+1} = (1-\alpha_i)\beta^i + \alpha_i \tbeta^i \ ,
\end{array}
\end{equation}

and notice that we recover the regular Frank-Wolfe update in the special case when $R(\cdot)$ is the indicator function $\bI_Q(\cdot)$ of a feasible region $Q$, see \cite{bach2015duality} and \cite{yu2017generalized} for a more detailed discussion on generalized Frank-Wolfe methods.  

\subsection{Related literature}

{\bf Stochastic Frank-Wolfe methods.}  There have been several lines of research that develop and investigate stochastic Frank-Wolfe methods, almost all of which are motivated by expected loss minimization in statistical and machine learning.  These methods by and large focus on the optimization problem:
	\begin{equation}
	\min_{\beta\in Q} f(\beta) = \tfrac{1}{n} \sum_{i=1}^n f_j(\beta) \ , 
	\end{equation}
where $Q$ is a closed and bounded convex set.  Here $f(\cdot)$ is the empirical risk, though some of the results we discuss below pertain to the more general ``infinite" setting of expected loss minimization where \begin{equation}\label{thc} f(\beta) = \int_\xi f(\beta;\xi) d\mu(\xi) \end{equation} for some appropriate probability measure $\mu(\cdot)$.

Hazan and Luo \cite{hazan2016variance} discuss several different stochastic Frank-Wolfe algorithms with increasing batch sizes over the course of iterations, including a straightforward stochastic Frank-Wolfe method (SFW), a stochastic variance reduced Frank-Wolfe method (SVRFW), and a stochastic variance-reduced conditional gradient sliding method (STORC).  In order to compute an $\varepsilon$-optimal solution, SFW requires $O(\tfrac{1}{\varepsilon^3})$ stochastic gradient calls (i.e., one computation of the gradient of one of the $f_j(\cdot)$ above) and $O(\tfrac{1}{\varepsilon})$ linear optimization oracle calls.  (Recall that we ignore other instance-specific constants such as sample size $n$ and/or variance measures of stochastic gradients, Lipschitz constants, curvature, and diameters constants, as our focus here is on the dependence on $\varepsilon$.)  SVRFW needs $O(\tfrac{1}{\varepsilon^2})$ stochastic gradient calls and $O(\tfrac{1}{\varepsilon})$ linear optimization oracle calls, as well as $O(\ln(\tfrac{1}{\varepsilon}))$ full gradient calls.  STORC (for STOchastic variance-Reduced Conditional gradient sliding) is a variance-reduced version of the conditional gradient sliding method of Lan and Zhou  \cite{lan2016conditional}.  The number of stochastic gradient calls that STORC needs is a function of some instance-specific properties: it is $O(\tfrac{1}{\varepsilon})$ if there is an optimal solution in the interior of $Q$, it is $O(\tfrac{1}{\varepsilon^{1.5}})$ more broadly, and is $O(\ln(\tfrac{1}{\varepsilon}))$ under strong convexity of $f(\cdot)$.  STORC also needs $O(\tfrac{1}{\varepsilon})$ linear optimization oracle calls as well as $O(\ln(\tfrac{1}{\varepsilon}))$ full gradient calls.

Lan and Zhou  \cite{lan2016conditional} present the stochastic conditional gradient sliding (SCGS) algorithm, which combines Nesterov's acceleration techniques and the Frank-Wolfe method.  In the absence of strong convexity, their stochastic Frank-Wolfe methodology requires $O(\tfrac{1}{\varepsilon^2})$ stochastic gradient calls and $O(\tfrac{1}{\varepsilon})$ linear optimization oracle calls.  And in the presence of strong convexity, their method requires $O(\tfrac{1}{\varepsilon})$ stochastic gradient calls and $O(\tfrac{1}{\varepsilon})$ linear optimization oracle calls.  

Mokhtari, Hassani and Karbasi  \cite{mokhtari2018stochastic} propose a stochastic Frank-Wolfe method called SCG (for stochastic conditional gradient) by introducing a momentum gradient estimator which does not require to increase the batch size.  SCG requires $O(\tfrac{1}{\varepsilon^3})$ stochastic gradient calls and $O(\tfrac{1}{\varepsilon^3})$ linear optimization oracle calls.  

In the online setting, Hazan and Kale \cite{hazan2012projection} proposed an online Frank-Wolfe method which requires $O(\frac{1}{\varepsilon^4})$ stochastic gradient calls and $O(\frac{1}{\varepsilon^4})$ linear optimization oracle calls.


Of the above mentioned methods, we point out as well that the algorithms SFW \cite{hazan2016variance}  as well as SCGS \cite{lan2016conditional} and SCG \cite{mokhtari2018stochastic} can also be implemented in the ``infinite'' setting of \eqref{thc}.

Unlike the deterministic Frank-Wolfe method, the above stochastic Frank-Wolfe methods do not achieve $O(\tfrac{1}{\varepsilon})$ gradient calls and $O(\tfrac{1}{\varepsilon})$ linear optimization oracle calls to achieve an $\varepsilon$-optimal solution, without some additional restrictive assumptions.  (This is the ``complexity gap'' mentioned in the abstract.)  From the above discussion, we see that such complexity is achieved only for STORC \cite{hazan2016variance} (when either (i) an optimal solution lies in the interior of $Q$, or (ii) when $f(\cdot)$ is strongly convex) and for SCGS \cite{lan2016conditional} when $f(\cdot)$ is strongly convex.  

\textcolor{black}{The primary motivation (and contribution) of this paper is to show that the typical setting of empirical loss minimization -- in which linear prediction is present -- is sufficient to guarantee that a suitable stochastic Frank-Wolfe method achieves both $O(\tfrac{1}{\varepsilon})$ gradient calls and $O(\tfrac{1}{\varepsilon})$ linear optimization oracle calls to achieve an $\varepsilon$-optimal solution.  Stated in the context of the literature on this topic, one does not need strong convexity nor interior optimal solutions to achieve the same complexity as deterministic Frank-Wolfe, so long as the setting is empirical loss minimization with linear prediction -- which is quite prevalent in such models.   At the same time, our proposed stochastic generalized Frank-Wolfe algorithm (Algorithm \ref{al:SFW}), which we call ``GSFW'' for short, has an additional factor of $n$ (the sample size) in the required number of stochastic gradient oracle calls and linear optimization oracle calls, namely $O(\tfrac{n}{\varepsilon})$, to compute an absolute $\varepsilon$-optimal solution of the empirical risk minimization problem with linear prediction.  This contrasts somewhat with the other algorithms in this suite (SCGM, SCGS, SFW) which do not have this extra factor of $n$.  In fact, GSFW does not necessarily dominate (nor is it dominated by) these other methods in computational complexity due to the differential appearance among this suite of methods of a variety of other constants (Lipschitz constants, curvature, diameter, stochastic gradient variance) in addition to $n$.  This is discussed further in Remark \ref{rem:factor_n}.}

\noindent {\bf Randomized Dual Coordinate Descent Methods.}  One of the interpretations of our stochastic Frank-Wolfe method is that it is a dual coordinate descent method in the dual space.  It thus is relevant to review the appropriate literature on dual coordinate descent in this context.  Dual coordinate descent methods have been widely used in statistical and machine learning applications. Stochastic dual coordinate ascent (SDCA) for solving \eqref{eq:primal} was first proposed in \cite{shalev2013stochastic}. There are many follow-up works on SDCA, for example, accelerated proximal randomized dual coordinate, see \cite{shalev2014accelerated}, \cite{lin2015accelerated}, using a non-uniform distribution to choose the coordinate \cite{qu2015quartz}, and a primal-dual coordinate method \cite{zhang2017stochastic}, among others. All of these dual methods (or primal-dual methods) require the regularizer $R(\cdot)$ to be a strongly convex function (or require adding a dummy strongly convex regularizer to the objective function). This contrasts with the standard Frank-Wolfe set-up where $R(\cdot)$ is an indicator function and so is not strongly convex.  Furthermore, in the Frank-Wolfe setup, the objective function in the dual problem \eqref{eq:dual} is not necessarily differentiable, which falls outside of the standard set-up for randomized coordinate descent \cite{richtarik2014iteration}, \cite{nesterov2012efficiency}.  (However, it turns out that as a byproduct of our analysis we obtain convergence guarantees for randomized coordinate descent applied to the (non-differentiable) dual problem \eqref{eq:dual}.) We further discuss the connections and differences between the above methods and our method in Appendix \ref{sec:connections-to-shalev}.

\noindent {\bf Variance Reduction Techniques for Stochastic Optimization.} There have been many recent developments of stochastic methods designed to directly tackle the optimization problem \eqref{eq:primal}. In order to obtain improved convergence guarantees over the standard Stochastic Gradient Descent (SGD) method, variance reduction techniques have been proposed and extensively studied in recent years. SAG \cite{schmidt2017minimizing} is the first variance reduction method in the literature that we are aware of. In contrast to the sublinear convergence rate of SGD, SAG and several concurrent and/or subsequent works -- such as SVRG \cite{johnson2013accelerating}, MISO \cite{mairal2015incremental}, and SAGA \cite{defazio2014saga} -- obtain linear convergence when the objective function is both smooth and strongly convex. Variance reduction techniques can also be applied to non-strongly convex optimization \cite{schmidt2017minimizing}, \cite{mairal2015incremental}, \cite{defazio2014saga}, \cite{allen2016improved}, which leads to improved convergence guarantees as well.  More recently, Allen-Zhu \cite{allen2017katyusha} has proposed an accelerated stochastic method for directly solving \eqref{eq:primal}.  We mention as well that the stochastic dual coordinate method \cite{shalev2013stochastic} also corresponds to a variant of a variance reduction technique in the primal space \cite{shalev2016sdca}.
We refer the reader to \cite{allen2017katyusha} for a more detailed discussion on variance reduction techniques overall.

\subsection{Contributions}

The overall contribution of this paper is Algorithm \ref{al:SFW}, which is a generalized stochastic Frank-Wolfe method (hence the moniker GSFW) designed to solve the empirical risk minimization problem with linear prediction \eqref{eq:primal}.   The specific contributions of Algorithm 1 and its analysis are as follows:

\begin{enumerate}

\item GSFW requires $O(\tfrac{1}{\varepsilon})$ stochastic gradient oracle calls and $O(\tfrac{1}{\varepsilon})$ linear optimization oracle calls to compute an absolute $\varepsilon$-optimal solution of the empirical risk minimization problem with linear prediction \eqref{eq:primal}, see Theorem \ref{thm:non-strong}. This in particular demonstrates that the typical setting of empirical loss minimization -- in which linear prediction is present -- is sufficient to guarantee that a suitable stochastic Frank-Wolfe method achieves the same complexity as its deterministic counterpart.  Stated in the context of the prevalent literature on stochastic Frank-Wolfe, one does not need strong convexity nor interior optimal solutions to achieve the same complexity as deterministic Frank-Wolfe, so long as the setting is empirical loss minimization with linear prediction -- which is quite prevalent in such models. \textcolor{black}{At the same time, GSFW has an additional factor of $n$ in the required number of stochastic gradient oracle calls and linear optimization oracle calls, namely $O(\tfrac{n}{\varepsilon})$, which contrasts (at least somewhat) with other comparable algorithms, see Remark \ref{rem:factor_n}.}

\item In the special case when $R(\cdot)$ is strongly convex, GSFW requires $O(\ln(\frac{1}{\varepsilon}))$ stochastic gradient oracle calls and $O(\ln(\frac{1}{\varepsilon}))$ linear optimization oracle calls to compute an absolute $\varepsilon$-optimal solution of \eqref{eq:primal}, see Theorem \ref{thm:strong}.

\item We show that GSFW is equivalent to a randomized coordinate mirror descent algorithm applied to the dual problem (Algorithm \ref{al:RCMD}), see Lemma \ref{lem:equi}.  Algorithm \ref{al:RCMD} can be viewed as a variant of the SDCA algorithm, and in this lens  our method can handle a non-strongly convex function $R(\cdot)$, in contrast to the current SDCA literature. 

\item Our work also implies a convergence bound for randomized coordinate mirror descent in the case when the objective function is the sum of a non-smooth function and a strongly convex separable function. This is discussed in Section \ref{sec:nonsmooth-cd}.

\item The recognition of the empirical risk minimization problem with linear prediction \eqref{eq:primal} as a problem of special interest due to the linear prediction structure.


\end{enumerate}

\subsection{Notation} We use $e_j$ to denote the $j^{\mathrm{th}}$ unit coordinate vector in $\RR^p$.  The $\ell_p$ norm is denoted $\|\cdot\|_p$.  We use $\dl_j(\cdot)$ to denote the first derivative of $l_j(\cdot)$.  The Bregman distance function associated with a convex function $h(\cdot)$ is defined as $\Dh(y,x):= h(y) - h(x) - \nabla h(x)^T(y-x)$.  We use $\EE$ to denote expectation and $\EE_{j_i}$ to denote expectation conditional on the randomly chosen index $j_i$.  For indicator functions, we use $\bI_Q(\cdot)$ to denote the indicator function for the set $Q$, namely $\bI_Q(\beta):=0$ if $\beta \in Q$, and $\bI_Q(\beta) :=+\infty$ otherwise; and we use $\bI_{\{\mathrm constraint\}}(\beta)$ to denote the indicator function of a particular constraint (or condition), namely $\bI_{\{\mathrm constraint\}}(\beta):=0$ if the constraint is true at $\beta$, and $\bI_{\{\mathrm constraint\}}(\beta) :=+\infty$ otherwise.  In a slight abuse of terminology we refer to the ``subgradient'' of a concave function when it is perhaps more technically accurate to refer to this as a sup-gradient.   A differentiable function $\ff$ is $\mu$-strongly convex with respect to a norm $\|\cdot\|$ if it holds that $f(y) \ge f(x) + \nabla f(x)^T(y-x) + \tfrac{\mu}{2}\|y-x\|^2$ for all $x,y \in \dom f(\cdot)$.  A differentiable function $\ff$ is $\mu$-strongly convex with respect to a reference function $h(\cdot)$ if it holds that $f(y) \ge f(x) + \nabla f(x)^T(y-x) + \mu \Dh(y,x)$ for all $x,y \in \dom f(\cdot)$.

\section{Dual problem, and equivalence of Algorithm \ref{al:SFW} in the dual with Randomized Coordinate Mirror Descent}\label{dualdual}
$$
f^*(y) := \sup_{x\in \dom f(\cdot)} \{y^T x - f(x)\} \ .
$$
We will also be interested in the following dual problem of \eqref{eq:primal} that is constructed using the conjugate functions of the component functions of  \eqref{eq:primal}:
\begin{equation}\label{eq:dual}
\mathrm{(D):} \ \ \ \ \ \ \ \max_{w} D(w):=-R^{*}\left(-\tfrac{1}{n}X^{T}w\right)-\tfrac{1}{n}\sum_{j=1}^{n}l_{j}^{*}(w_{j}) \ .
\end{equation}
Notice that we can write:
\begin{equation}\label{eq:altered}R^{*}\left(-\tfrac{1}{n}X^{T}w\right)= - \min_{\beta}\left\{\tfrac{1}{n} w^{T}X\beta+R(\beta)\right\}  \ .
\end{equation}
Also, defining the convex/concave saddle-function $\phi(\cdot,\cdot)$:
\begin{equation}\label{eq:primal-dual}
\phi(\beta,w):=\tfrac{1}{n} w^{T}X\beta-\tfrac{1}{n}\sum_{i=1}^{n}l_{i}^{*}(w_{i})+R(\beta) \ ,
\end{equation} we can write (P) and (D) in saddlepoint minimax format as:

\begin{equation}\label{eq:minmax} \mathrm{(P):} \ \ \min_{\beta} \max_{w} \phi(\beta,w) \ \ \ \ \  \ \ \mathrm{and} \ \ \ \ \ \ \ \  \mathrm{(D):} \ \ \max_{w}\min_{\beta} \phi(\beta,w) \ .
\end{equation}

Another standard first-order method for convex optimization is the mirror descent algorithm (also called primal gradient method with Bregman distance) \cite{tsengaccelerated}, \cite{lu2018relatively}, \cite{lu2017relative}, \cite{beckteb03mirror}, which we now briefly review in the context of solving the dual problem (D) in \eqref{eq:dual}, which is a concave maximization problem.  The Bregman distance of a differentiable ``prox'' function $h(\cdot)$ is defined to be:
$$D_h(w_1,w_2):=h(w_1)-h(w_2)-\langle \nabla h(w_2), w_1 - w_2 \rangle \ . $$ The (deterministic) mirror descent algorithm for solving (D) has the following update:
$$
w^{i+1} \gets \arg\min_{w} \{ -\eta_i g(w^i)^T (w-w^i) +  \Dh(w,w^i) \} \ ,
$$
where $g(\cdot)$ is a subgradient of the objective function $D(\cdot)$ at $w$ (which we call a subgradient even though $D(\cdot)$ is concave), and $\{\eta_i\}$ is the step-size sequence.  It is shown in Bach \cite{bach2015duality} that the generalized Frank-Wolfe method for the primal \eqref{eq:primal} is equivalent to mirror descent algorithm for the dual \eqref{eq:dual}.

Algorithm \ref{al:RCMD} presents a Randomized Coordinate Mirror Descent method applied to solve the dual problem D. The algorithm uses the average of the conjugate functions $l_i^*(\cdot)$ as the prox function, namely $h(\cdot)=\frac{1}{n}\sum_{i=1}^{n}l_{i}^{*}(w_{i})$, and it initializes the dual variable $w^0$ to be the prox-center (which is the point that minimizes the prox function). At the start of the $i^{\mathrm{th}}$ iteration, the algorithm randomly chooses a coordinate $j_i$ and computes the $j_i^{\mathrm{th}}$ coordinate of a subgradient of the dual objective function $D(w)$ at $w=w^i$, since indeed it is straightforward to verify that $\frac{1}{n}(X\tbeta^{i}-\nabla L^{*}(w^{i}))$ is a subgradient of $D(w)$ at $w=w^i$.  The algorithm then performs a coordinate mirror descent step to update the dual variable $w^i$.  Last of all -- and optionally since it does not affect future computations -- the algorithm updates the primal variable $\bbeta^i$ in order to compute a primal-dual optimality gap certificate.\medskip

\begin{algorithm}
\caption{Randomized Coordinate Mirror Descent applied to the dual problem \eqref{eq:dual}}\label{al:RCMD}
$ \ $
\begin{algorithmic}
\STATE {\bf Initialize.}  Define the prox function $h(w):=\frac{1}{n}\sum_{i=1}^{n}l_{i}^{*}(w_{i})$.  Initialize with $w^{0}=\arg\min_{w} \tfrac{1}{n}\sum_{i=1}^{n}l_{i}^{*}(w_{i})$ and
step-size sequences $\{\alpha_{i}\} \in (0,1]$ and $\{\eta_{i}\} \in (0,1]$. (Optional: set $\bbeta^{-1}=0$.)\medskip

\STATE  For iterations $i=0,1, \ldots $

\ \ \ \ \ \ \  {\bf Compute Randomized Coordinate of Subgradient of $D(\cdot)$ at $w^i$}\\

\ \ \ \ \ \ \ \ \ \ \ \ Compute $\tbeta^{i} \in \arg\min_{\beta}\left\{\left( \tfrac{1}{n}(w^i)^T X \beta+R(\beta) \right) \right\}$ \\ \medskip

\ \ \ \ \ \ \ \ \ \ \ \    {\bf Choose random index.}  Choose $j_{i} \in {\cal U}[1, \ldots, n]$ \\
\ \ \ \ \ \ \ \ \ \ \ \  {\bf Compute subgradient coordinate vector:} $\tilde{g}^i \gets \frac{1}{n}\left(x_{j_{i}}^{T}\tbeta^{i}-\dot{l}_{j_{i}}^{*}(w_{j_{i}}^{i})\right)e_{j_{i}}$ \\ \medskip

\ \ \ \ \ \ \    {\bf Update dual variable:} Compute $w^{i+1}=\arg\min_w \left\{ \left\langle -\eta_{i}\tilde g^i,w-w^{i}\right\rangle +D_{h}(w,w^{i})\right\} $ \\ \medskip

\ \ \ \ \ \ \  {\bf (Optional Accounting:)} $\bbeta^{i} \leftarrow (1-\alpha_i) \bbeta^{i-1} + \alpha_i \tbeta^i$. \medskip

\end{algorithmic}
\end{algorithm}\medskip

The main result of this section is the following lemma concerning the equivalence of Algorithm \ref{al:SFW} and Algorithm \ref{al:RCMD}.\medskip

\begin{lem}\label{lem:equi} {\bf (Equivalence Lemma)}  Algorithm \ref{al:SFW} and Algorithm \ref{al:RCMD} are equivalent as follows:  the iterate sequence of either algorithm exactly corresponds to an iterate sequences of the other. \qed
\end{lem}\medskip\medskip


As a means to proving the lemma, we first reinterpret the update of $w_{j_i}$ at iteration $i$ of Algorithm \ref{al:RCMD} in the following proposition:\medskip

\begin{prop}\label{prop:wupdate} At iteration $i$ of Algorithm \ref{al:RCMD} it holds that:

(1.) $\dlsji(w_{j_{i}}^{i+1})=(1-\eta_{i})\dot{l}_{j_{i}}^{*}(w_{j_{i}}^{i}) + \eta_{i}x_{j_{i}}^{T}\tbeta^{i}$, and

(2.) $w_{j_{i}}^{i+1}=\dot{l}_{j_{i}}\left((1-\eta_{i})\dot{l}_{j_{i}}^{*}(w_{j_{i}}^{i})+ \eta_{i}x_{j_{i}}^{T}\tbeta^{i}\right) $.
\end{prop}
{\bf Proof:}
Because $h(w)$ is a coordinate-wise separable function, we can rewrite the update for $w_{j_{i}}^{i+1}$ as
\begin{equation*}
\begin{array}{ll}
w_{j_{i}}^{i+1} & =\arg\min_{w_{j_{i}}} \left\langle -\frac{\eta_{i}}{n}\left(x_{j_{i}}^{T}\tbeta^{i}-\dot{l}_{j_{i}}^{*}(w_{j_{i}}^{i})\right),w_{j_{i}}\right\rangle +D_{\frac{1}{n} l_{j_{i}}^{*}}(w_{j_{i}},w_{j_{i}}^{i})\\ \\

& =\arg\min_{w_{j_{i}}}\left\langle -\eta_{i}\left(x_{j_{i}}^{T}\tbeta^{i}-\dot{l}_{j_{i}}^{*}(w_{j_{i}}^{i})\right),w_{j_{i}}\right\rangle +D_{ l_{j_{i}}^{*}}(w_{j_{i}},w_{j_{i}}^{i})\\ \\

& =\arg\min_{w_{j_{i}}}\left\langle -\eta_{i}x_{j_{i}}^{T}\tbeta^{i}-(1-\eta_{i})\dot{l}_{j_{i}}^{*}(w_{j_{i}}^{i}),w_{j_{i}}\right\rangle +l_{j_{i}}^{*}(w_{j_{i}}) \ .  \end{array}
 \end{equation*}

From the first-order optimality condition of the above $1$-dimensional problem we have $\dls_{j_i}(w_{j_i}^{i+1})=\eta_{i}x_{j_{i}}^{T}\tbeta^{i}+(1-\eta_{i})\dot{l}_{j_{i}}^{*}(w_{j_{i}}^{i})$, which shows {\it (1.)}; and {\it (2.)} follows directly from {\it (1.)} by the properties of the conjugate function in Proposition \ref{prop:gradient_conjugate}. \qed

%
%
%


{\bf Proof of Lemma \ref{lem:equi}}  We show that the iterate sequence of Algorithm \ref{al:RCMD} corresponds exactly to an iterate sequence of Algorithm \ref{al:SFW}.  The $\{s^i\}$ sequence is not formally defined in Algorithm \ref{al:RCMD}, so let us define $s^i := \nabla L^*(w^i)$ for all $i = 0, 1, \ldots$, which is consistent through conjugacy with the relationship $w^i = \nabla L(s^i)$ in the Optional Accounting step of Algorithm \ref{al:SFW} (see Proposition \ref{prop:gradient_conjugate}).  In order to show the correspondence we proceed by induction on the iteration counter $i$.  For $i=0$ we have from conjugacy that $s^0 := \nabla L^*(w^0) = 0 $ from the definition $w^0$ in the  initialization of Algorithm \ref{al:RCMD}.  We also need to show that $\tbeta^0$ is a solution to the linear optimization oracle problem in Algorithm \ref{al:SFW}.  We have for all $i=0, \ldots$, that:
$$\begin{array}{rcl}
\tbeta^{i} \in \arg\min_\beta\left\{ \tfrac{1}{n}\left(w^{i}\right)^{T}X\beta+R(\beta)\right\}    &=& \arg\min_\beta\left\{ \tfrac{1}{n}\left(\nabla L(s^{i})\right)^{T}X\beta+R(\beta)\right\} \\ \\
&=& \arg\min_\beta\left\{ \left(d^{i}\right)^{T}\beta+R(\beta)\right\} \ ,
\end{array}$$thus showing that $\beta^{i}$ corresponds to a linear optimization oracle solution at iteration $i$ in Algorithm \ref{al:SFW} for all $i=0, \ldots$.  Now suppose that the correspondence holds for some iteration counter $i$, and let us examine $s^{i+1}:= \nabla L^*(w^{i+1})$.  We have from Proposition \ref{prop:wupdate} that:
\begin{equation}\label{robo3}
w_{j_{i}}^{i+1} =\dot{l}_{j_{i}}\left((1-\eta_{i})\dot{l}_{j_{i}}^{*}(w_{j_{i}}^{i}) + \eta_{i}x_{j_{i}}^{T}\tbeta^{i} \right)=\dot{l}_{j_{i}}\left((1-\eta_{i})s_{j_{i}}^{i} + \eta_{i}x_{j_{i}}^{T}\tbeta^{i} \right) \ ,
\end{equation}
where the first equality is from Proposition \ref{prop:wupdate} and the second equality uses induction.  This then implies that $$ (1-\eta_{i})s_{j_{i}}^{i} + \eta_{i}x_{j_{i}}^{T}\tbeta^{i} = \dl (w_{j_{i}}^{i+1}) = s^{i+1} \ . $$
And for all coefficient indices $j \ne i $ we have
$$s_j^{i+1} = \dl^*(w^{i+1}_j) = \dl(w^i_j)= s^i \ , $$
where the second equality follows from conjugacy, whereby $s^{i+1}$ satisfies the update rule as stated in Algorithm \ref{al:SFW}, thus demonstrating that the iterate sequence of Algorithm \ref{al:RCMD} corresponds exactly to an iterate sequence of Algorithm \ref{al:SFW}.  The same type of analysis as above can be used to prove that the iterate sequence of Algorithm \ref{al:SFW} corresponds exactly to an iterate sequence of Algorithm \ref{al:RCMD}. \qed

\section{Convergence Guarantees}\label{convergenceguarantee}
In this section we develop computational guarantees for Algorithm \ref{al:RCMD}, which automatically provide computational guarantees for Algorithm \ref{al:SFW} due to the equivalence shown in Theorem \ref{lem:equi}.  Our first -- and main -- result is Theorem \ref{thm:non-strong}, which is an expected $O(1/k)$ guaranteed decrease in the duality gap between (P) and (D).  Secondly, in the case when $R(\cdot)$ is a strongly convex function, we present a linear convergence result on the duality gap in Theorem \ref{thm:strong}.  We start by defining two measures -- $M$ and $D_{\max}$ -- associated with (P) and whose values will enter our computational bounds.

Let $M := \max_{\beta \in \dom R(\cdot)} \max_{j=1, \ldots, n} \{| x_j^T\beta | \}$, and note that $M < +\infty$ since $\dom R(\cdot)$ is bounded by Assumption \ref{ass:smooth}.\medskip

Let ${\cal W} \subset \RR^n$ be the set of ``optimal $w$ responses'' to values $\beta \in \dom R(\cdot)$ in the saddle-function $\phi(\beta, w)$, namely: $${\cal W} := \{ \hw \in \RR^n : \hw \in \arg \max_w \phi(\hbeta, w) \ \mbox{for~some~} \hbeta \in \dom R(\cdot) \} \ , $$

and let $D_{\max}$ be any upper bound on $D_h (\hw,w^0)$ as $\hw$ ranges over all values in ${\cal W}$, so that
$$
D_h (\hw ,w^0)\le D_{\max} \ \  \mbox{for~all~}  \hw \in \cal W \ .
$$
Note at the moment that there is no guarantee that $D_{\max} < +\infty$, but this will be remedied below in Proposition \ref{prop:Dmax}.\medskip

\begin{prop}\label{prop:Dmax}
	Under Assumption \ref{ass:smooth} it holds that $D_{\max} \le \gamma M^2$.
\end{prop}

Before proving this proposition, we first show that there is a natural boundedness constraint for the dual problem:\medskip

\begin{prop}\label{prop:well-behavior-of-w}
	Let $T:=\left\{ w \in \RR^n : |\|w-w^{0}\|_{\infty}\le \gamma M\right\}$.  Then:
	\begin{enumerate}
		\item for any $\hbeta\in \dom R(\cdot)$ it holds that $\argmax_w \phi (\hbeta, w) \in T$, and
		\item for all $w^i$ generated in Algorithm \ref{al:RCMD}, it holds that  $w^i \in T$.
	\end{enumerate}
\end{prop}

\textbf{Proof.} We first prove {\it (1.)}.  Notice that $w_{0}= \nabla L(0)$ and $\argmax_w \phi (\hbeta, w) = \nabla L(X\hbeta)$ (from conjugacy via Proposition \ref{prop:gradient_conjugate}) , whereby the $\gamma$-smoothness of $l_{j}(\cdot)$ implies that
$$
\left\|\arg\max_{w}\phi(\hat{\beta},w)-w_{0}\right\|_{\infty}=\left\| \nabla L(X\hbeta)-\nabla L(0)\right\|_{\infty}=\max_j \left|\dl_j(x_j^T\hat \beta)-\dl_j(0)\right|\le\gamma\max_{j}|x_{j}^{T}\hat \beta| \le \gamma M \ ,
$$

which proves {\it (1.)}. It follows from Proposition \ref{cor:bound} that for any coordinate $j$ and iterate $i$ it holds that $\left|\dot{l}_{j}^{*}(w_{j}^{i})\right|\le M$. Together with $\dls_j(w_j^0)=0$, we have

$$
\tfrac{1}{\gamma}\left|w_{j}^{i}-w_{j}^{0}\right|\le \left|\dot{l}_{j}^{*}(w_{j}^{i})-\dot{l}_{j}^{*}(w_{j}^{0}) \right|\le M \ ,
$$
(where the first inequality is from the $\frac{1}{\gamma}$-strong convexity of $l_j^*(w_j)$), from which it follows that $\|w^i - w^0\|_\infty \le \gamma M$, which proves {\it (2.)}. \qed

{\bf Proof of Proposition \ref{prop:Dmax}:}  Let $L(s):=\sum_{j=1}^n l_j(s_j)$ and $L^*(w):=\sum_{j=1}^n l^*_j(w_j)$, and note that $L(\cdot)$ and $L^*(\cdot)$ are a conjugate pair.  Let $\hw \in \cal W$ and let $\hbeta$ be such that $\hw \in \arg\max_w \phi(\hbeta, w)$.  Then


\begin{equation*}
\begin{array}{lcl}
D_{h}(\hw, w^0) & = &\tfrac{1}{n}\left( L^*(\hw)-L^*(w^0)\right)\\ \\
& =&\tfrac{1}{n}\left( \left(\hw \right)^{T}X\hbeta-L(X\hbeta) -L^{*}(w^{0})\right)\\ \\
& \le&\tfrac{1}{n}\left(\max_{w\in T,\beta\in \dom R(\cdot)}\left\{ w^{T}X\beta-L(X\beta)\right\} -L^{*}(w^{0})\right)\\ \\
& =&\tfrac{1}{n}\left(\max_{w\in T,\beta\in  \dom R(\cdot)}\left\{ (w-w^{0})^{T}X\beta+\left(w^{0}\right)^{T}X\beta-L(X\beta)\right\} -L^{*}(w^{0})\right)\\ \\
& \le&\tfrac{1}{n}\left(\max_{w\in T,\beta\in  \dom R(\cdot)}\left\{ (w-w^{0})^{T}X\beta\right\} +\max_{\beta\in  \dom R(\cdot)}\left\{ \left(w^{0}\right)^{T}X\beta-L(X\beta)\right\} -L^{*}(w^{0})\right)\\ \\
& \le &\tfrac{1}{n}\left(n\max_{w\in T,\beta\in  \dom R(\cdot)}\|w-w^{0}\|_{\infty}\|X\beta\|_{\infty}+L^{*}(w^{0})-L^{*}(w^{0})\right)\\ \\
& \le &\gamma M^{2} \ ,
\end{array}
\end{equation*}

where the second equality follows from Proposition \ref{prop:gradient_conjugate}, the first inequality uses $\hbeta \in  \dom R(\cdot)$ and $\hw \in T$ (from Proposition \ref{prop:well-behavior-of-w}), and the last inequality uses $\max_{\beta\in  \dom R(\cdot)}\|X\beta\|_{\infty}\le M$. \qed \medskip

\begin{rem}\label{remarkrob}  A suitable value of $D_{\max}$ can often be easily derived based on the structure of $l_j(\cdot)$.  For example, in logistic regression where the loss function is $l_j(s_j):=\log(1+\exp(-y_j s_j))$  for the given label $y_j \in \{-1,1\}$, we have $l_j^*(w_j)=-y_j w_j \ln(-y_jw_j) + (1+y_j w_j)\ln(1+y_j w_j)$ with $\dom l_j^*(\cdot) = \{w_j : 0 \le -y_jw_j \le 1 \}$ (where $a\ln(a) :=0$ for $a=0$). Therefore for all $\hw \in \cal W$ it holds that $$D_h (\hw,w^0)\le \max_{0\le -Yw\le e} D_h (w,w^0)=\tfrac{1}{n}\left(\max_{0\le -Yw\le e}L^*(w)-L^*(w^0)\right)=\ln (2) \ ,$$
where $Y$ is the diagonal matrix whose diagonal coefficients correspond to $y$ and $e=[1,\ldots,1]^T$, so we may set $D_{\max}=\ln (2)$.
\end{rem}

Notice in Algorithm \ref{al:SFW} and Algorithm \ref{al:RCMD} that $j_i$ is a random variable; and that $s^i$, $d^i$, $w^i$, etc., are random variables that depend on all previous random variable values $j_0, j_1, \ldots, j_{i-1}$, and we denote this string of random variables by 
\begin{equation}\label{eq:xi}
\xi_i=\{j_0, j_1, \ldots, j_{i-1}\} \ .
\end{equation}

We now state our main computational guarantee for Algorithm \ref{al:RCMD} (and hence for Algorithm \ref{al:SFW} as well).\medskip

\begin{thm}\label{thm:non-strong}
Consider the Stochastic Generalized Frank-Wolfe method (Algorithm \ref{al:SFW}) or the Randomized Dual Coordinate Mirror Descent method (Algorithm \ref{al:RCMD}), with step-size sequences $\alpha_i=\frac{2(2n+i)}{(i+1)(4n+i)}$ and $\eta_i=\frac{2n}{2n+i+1}$ for $i=0, 1, \ldots$.  Denote $$\bw^{k} = \frac{2}{(4n+k)(k+1)}\sum_{i=0}^k(2n+i) w^i \ . $$
Under Assumption \ref{ass:smooth}, it holds for all $k \ge 0$ that
$$
\EE_{\xi_k}\left[ P(\bbeta^k) - D(\bw^k)\right] \le \frac{8n\gamma M^{2}}{\left(4n+k\right)}+\frac{2n(2n-1)D_{\max}}{(4n+k)(k+1)}  \ \le \ \frac{8n\gamma M^{2}}{\left(4n+k\right)}+\frac{2n(2n-1)\gamma M^2}{(4n+k)(k+1)} \ . \ \qed
$$
\end{thm}\medskip

\begin{rem}\label{rem:factor_n} \textcolor{black}{If we include the dependence on the sample size $n$ in the big-O notation, then Theorem \ref{thm:non-strong} shows that GSFW requires $O(\tfrac{n}{\varepsilon})$ iterations to compute an absolute $\varepsilon$-optimal solution of the empirical risk minimization problem with linear prediction \eqref{eq:primal}, and indeed that is the same order of gradient computations (over individual samples) as the deterministic Frank-Wolfe method \cite{jaggi2013revisiting}. Remark \ref{wilderness} discusses how this dependency changes in the presence of mini-batches, and shows how the convergence rate changes as the algorithm morphs from stochastic to deterministic as the mini-batch size is increased.  Actually the results herein have a similar structure as randomized coordinate descent (RCD) complexity bounds for solving smooth optimization (see, e.g, Theorem 5 in \cite{nesterov2012efficiency}), wherein RCD requires $O(\tfrac{n}{\varepsilon})$ to obtain an $\varepsilon$-optimal solution and $n$ therein refers to the number of coordinates.  As we can see in Theorem \ref{thm:non-strong}, GSFW has a superior convergence rate in $\varepsilon$ compared to the existing stochastic Frank-Wolfe methods (such as SFW, SCGM, SCGS and SVRF), but has the additional factor of $n$ in contrast with these other methods.  In fact, GSFW does not necessarily dominate (nor is it dominated by) these other methods in computational complexity due to the differential appearance among this suite of methods of a variety of other constants (Lipschitz constants, curvature, diameter, stochastic gradient variance) in addition to $n$.}
\end{rem}

The following string of propositions will be needed for the proof of Theorem \ref{thm:non-strong}.\medskip

\begin{prop}
For all iterates $i$ and any $j\in\left\{ 1,\ldots,n\right\} $ it holds that $\left|\dot{l}_{j}^{*}(w_{j}^{i})\right|\le M$.
\end{prop}

{\bf Proof.} We prove this by induction on $i$. The proposition is true for $i=0$ because $\dot{l}_{j}^{*}(w_{j}^{0})=0$ for all $j$ by the definition of $w^0$. Next suppose that $\left|\dot{l}_{j}^{*}(w_{j}^{i})\right|\le M$ for a given iterate $i$ and for all $j=1, \ldots, n$. Then at iteration $i+1$ and any $j\not=j_{i}$ we have  $w_j^{i+1}=w_j^i$, whereby $\left|\dot{l}_{j}^{*}(w_{j}^{i+1})\right|=\left|\dot{l}_{j}^{*}(w_{j}^{i})\right|\le M$.  And it follows from Proposition \ref{prop:wupdate} that
\[
\left|\dot{l}_{j_{i}}^{*}(w_{j_{i}}^{i+1})\right|=\left|\left(1-\eta_{i}\right)\dot{l}_{j_{i}}^{*}(w_{j_{i}}^{i})+\eta_{i}x_{j_{i}}^T\tbeta^{i}\right|\le(1-\eta_{i})M+\eta_{i}M=M \ ,
\]
and therefore for any $j=1, \ldots, n$, we have $\left|\dot{l}_{j}^{*}(w_{j}^{i+1})\right|\le M$, which completes the proof by induction. \qed

As a simple corollary we obtain an upper bound on $\|\tilde{g}^i\|_{2}$ as follows:\medskip

\begin{cor}\label{cor:bound}
$\|\tilde{g}^i\|_{2} =\frac{1}{n}\left| x_{j_{i}}^T \tbeta^i - \dlsji(\wiji) \right| \le \frac{2M}{n}$.
\end{cor}\medskip

\begin{prop}\label{prop:strong}
$h(\cdot)$ is
$\frac{1}{n\gamma}$-strongly convex with respect to the norm $\|\cdot\|_{2}$.
\end{prop}
{\bf Proof.} Recall that $h(w)=\frac{1}{n}\sum_{j=1}^{n}l_{j}^{*}(w_{j})$. It follows from Assumption \ref{ass:smooth} and Proposition \ref{prop:gradient_conjugate} that $\dls_j(\cdot)$ is $\frac{1}{\gamma}$-strongly convex. Therefore for any $w^1, w^2 \in \dom h(\cdot)$ it holds that:
\begin{equation*}
\begin{array}{lcl}
h(w^1)& = &\frac{1}{n}\sum_{j=1}^{n}l_{j}^{*}(w^1_{j}) \\ \\
& \ge & \frac{1}{n}\sum_{j=1}^{n}\left( l_{j}^{*}(w^2_{j}) + \dls_j(w_j^2)(w_j^1-w_j^2)+\frac{1}{2\gamma}|w^2_j-w^1_j|^2\right) \\ \\
& = & h(w^2) + \langle\nabla h(w^2), w^1-w^2\rangle + \frac{1}{2n\gamma}\|w^2-w^1\|_2^2 \ . \ \qed
\end{array}
\end{equation*}

\medskip

\begin{prop}\label{prop:equality}
$\phi(\tbeta^{i},w)=D(w^{i})+\left\langle \nabla_{w}\phi(\tbeta^{i},w^{i}),w-w^{i}\right\rangle -D_{h}(w,w^{i}).$
\end{prop}
\textbf{Proof. } The proof follows from straightforward substitution using $\phi(\tbeta^{i},w)=\frac{1}{n} \left(w^{T}X\tbeta^{i}-\sum_{j=1}^n l_{j}^{*}(w_{j})\right) + R(\tbeta^i)$ and noticing from the construction of $\tbeta^i$ that $D(w^i)=\phi(\tbeta^i, w^i)$. \qed\medskip

\begin{prop}\label{robfriday}  Consider the series $\{ \alpha_i \}$ defined by $\alpha_i = \frac{2(2n+i)}{(4n+i)(i+1)}$ for $i \ge 0$ and define the series $\{\bbeta^i\}$ by $\bbeta^{-1} =0$ and $\bbeta^i = (1-\alpha_i)\bbeta^{i-1} + \alpha_i \tbeta^i$ for $i\ge 0$.  Also define $\gamma_i = 2n+i$ for $i \ge 0$.  Then
$$ \bbeta^k = \frac{\sum_{i=0}^k \gamma_i \tilde{\beta}^i}{\sum_{i=0}^k \gamma_i} \ \ \mbox{for~all~} k \ge 0 \ . $$
\end{prop}
\noindent {\bf Proof:} The proof follows easily by induction and using $\sum_{i=0}^k \gamma_i = \frac{(4n+k)(k+1)}{2}$. \qed

{\bf Proof of Theorem \ref{thm:non-strong}.}  Denote $g^i:=\frac{1}{n} \left( X\tbeta^i-\nabla L^*(w^i)\right)$, whereby $g^i$ is a subgradient of $D(w)$ at $w^i$, and $\tg^i$ is an unbiased estimator of $g^i$ up to the scalar $n$, namely $\EE_{j_i} [ \tg^i ]=\frac{1}{n}g^i$.  Therefore we have for any $i$ and any $w \in \cal W$ that:
\begin{equation}\label{eq:nonsmooth-theonehand}
\begin{array}{ll}
\left\langle -g^i,w-w^{i}\right\rangle  & =n \EE_{j_i} \left[\left\langle -\tg^i,w-w^{i}\right\rangle\right] \\ \\
 & \ge n \EE_{j_i} \left[\left\langle -\tg^i,w^{i+1}-w^{i}\right\rangle +\frac{1}{\eta_{i}}D_{h}(w^{i+1},w^{i})+\frac{1}{\eta_{i}}D_{h}(w,w^{i+1})-\frac{1}{\eta_{i}}D_{h}(w,w^{i})\right]\\ \\
 & \ge n \EE_{j_i} \left[\left\langle -\tg^i,w^{i+1}-w^{i}\right\rangle +\frac{1}{2n\gamma\eta_{i}}\|w^{i+1}-w^{i}\|_{2}^{2}+\frac{1}{\eta_{i}}D_{h}(w,w^{i+1})-\frac{1}{\eta_{i}}D_{h}(w,w^{i})\right]\\ \\
  & \ge n \EE_{j_i} \left[ -\frac{1}{2}n\gamma \eta_{i} \|\tg^i\|_2^2+\frac{1}{\eta_{i}}D_{h}(w,w^{i+1})-\frac{1}{\eta_{i}}D_{h}(w,w^{i})\right]\\ \\
 & \ge -2\gamma M^{2}\eta_{i}+\frac{n}{\eta_{i}}\EE_{j_i} [D_{h}(w,w^{i+1})]-\frac{n}{\eta_{i}}D_{h}(w,w^{i}) \ ,
\end{array}
\end{equation}
where the first inequality is from the ``three point property'' of Tseng (Lemma \ref{tpp} in the Appendix), the second
inequality is due to the fact that $h(w)$ is $\frac{1}{n\gamma}$-strongly convex with respect to the norm $\|\cdot\|_2$ (Proposition \ref{prop:strong}), and the third inequality is an application of the basic inequality $\langle x, y\rangle \le \tfrac{1}{2}\|x\|_2^2 +  \tfrac{1}{2}\|y\|_2^2$,  and the last inequality uses Corollary \ref{cor:bound}.

On the other hand, we have from Proposition \ref{prop:equality} that
\begin{equation}\label{eq:nonsmooth-theotherhand}
\left\langle -g^i,w-w^{i}\right\rangle   =\left\langle -\nabla_{w}\phi(\tbeta^{i},w^{i}),w-w^{i}\right\rangle  =D(w^{i})-\phi(\tbeta^{i},w)-D_{h}(w,w^{i}) \ .
\end{equation}

Combining \eqref{eq:nonsmooth-theonehand} and \eqref{eq:nonsmooth-theotherhand} and rearranging yields
$$
-(\phi(\tbeta^{i},w)-D(w^{i}))\ge-2\gamma M^{2}\eta_{i}+\tfrac{n}{\eta_{i}}\EE_{j_i}[D_{h}(w,w^{i+1})]-\left(\tfrac{n}{\eta_{i}}-1\right)D_{h}(w,w^{i}).
$$
Substituting $\eta_{i}=\frac{2n}{2n+i+1}$ and multiplying by $2n+i$ results, we arrive at the following inequality after rearranging terms:
$$\begin{array}{l}
(2n+i)(\phi(\tbeta^{i},w)-D(w^{i}))\\
\le  \ \ 4n\gamma  M^{2}\left(\tfrac{2n+i}{2n+i+1}\right)+\tfrac{1}{2}\left((2n+i)(2n+i-1)D_{h}(w,w^{i})-(2n+i)(2n+i+1)\EE_{j_i}[ D_{h}(w,w^{i+1})] \right) \ .
\end{array}$$
Summing the above inequality for $i=0,\ldots,k$ and recalling from Proposition \ref{robfriday} that $\bar{\beta}^{k}:=\frac{2}{(4n+k)(k+1)}\sum_{i=0}^{k}(2n+i)\tbeta^{i}$,
and taking the unconditional expectation over $\xi_k$ (recall the definition of $\xi_k$ in \eqref{eq:xi}), we arrive at:
$$\begin{array}{rcl}
\frac{(4n+k)(k+1)}{2}\EE_{\xi_k}[\phi(\bar{\beta}^{k},w)-D(\bar{w}^{k})] & = & \left(\sum_{i=0}^k 2n+i \right)\EE_{\xi_k}[\phi(\bar{\beta}^{k},w)-D(\bar{w}^{k})] \\ \\
& \le & \EE_{\xi_k}\left[\sum_{i=0}^{k}(2n+i)(\phi(\tbeta^{i},w)-D(w^{i}))\right]\\ \\
& \le & 4(k+1)n\gamma M^{2}+\tfrac{1}{2}(2n)(2n-1)D_{h}(w,w^{0})\\ \\
& \le & 4(k+1)n\gamma M^{2}+n(2n-1)D_{\max} \ ,
\end{array}$$
where the first inequality uses the convexity of $\phi(\beta,w)$
over $\beta$ and the concavity of $D(w)$, the second inequality
follows from the summation and canceling terms in the telescoping series, and the third inequality uses $w \in \cal W$.  Choosing $\hw=\arg\max_{w}\phi(\bar{\beta}^{k},w)$,
we have $P(\bar{\beta}^{k})=\phi(\bar{\beta}^{k},\hw)$, which yields:
\[
\EE_{\xi_k}[P(\bar{\beta}^{k})-D(\bar{w}^{k})]\le\frac{8n\gamma M^{2}}{\left(4n+k\right)}+\frac{2n(2n-1)D_{\max}}{(4n+k)(k+1)} \ ,
\]thus showing the first inequality in the statement of the theorem.
The second inequality in the statement of the theorem then follows as a simple application of Proposition \ref{prop:Dmax}. \qed

\subsection{Linear Convergence when $R(\cdot)$ is Strongly Convex}\label{sec:linear_conv}

In this section, we further assume $R(\cdot)$ is a $\mu$-strongly convex function, and we develop a linear convergence guarantee for Algorithms \ref{al:SFW} and \ref{al:RCMD} .  We first formally define a separable function.\medskip

\begin{mydef}
The function $h:R^{n}\rightarrow R$ is separable if
\[
h\left(x\right)=\sum_{i=1}^{n}h_{i}\left(x_{i}\right),
\]
where $x_{i}$ is the $i^{\mathrm{th}}$ coordinate of $x$ and $h_{i}$ is a univariate function.
\end{mydef}\medskip

Next we introduce the notation of relative smoothness and relative strong convexity developed recently in \cite{lu2018relatively}\cite{lu2017relative}\cite{bauschke2016descent}\cite{hanzely2018fastest}\cite{van2017forward}. We adapt a simplified version of the coordinate-wise relative smoothness condition as in \cite{hanzely2018fastest}.\medskip

\begin{mydef}
$f(\cdot)$ is coordinate-wise $\sigma$-smooth relative to a separable reference function $h(\cdot)$ if for any $x$, scalar $t$ and coordinate $j$ it holds that:\begin{equation}\label{eq:coor-smoothness}
f(x+t e_j)\le f(x) + \langle \nabla f(x), t e_j \rangle + \sigma D_h(x+t e_j,x) \ .
\end{equation}
\end{mydef}

We also adapt the notion of relative strong convexity developed in \cite{lu2018relatively}.\medskip

\begin{mydef}
$f(\cdot)$ is $\mu$-strongly convex relative to $h(\cdot)$ if for any $x,y$, it holds that
\begin{equation}
f(y)\ge f(x) + \langle \nabla f(x), y-x \rangle + \mu D_h(y,x) \ .
\end{equation}
\end{mydef}\medskip

The next proposition states that the dual function $D(w)$ is both coordinate-wise smooth and strongly concave relative to the reference function $h(w):=\frac{1}{n}\sum_{j=1}^n \ls_j(w_j)$.  In the proposition, recall that $x_j$ is the $j^{\mathrm{th}}$ row of the matrix $X$.\medskip

\begin{prop}\label{prop:coordinate smoothness}
$  $

(1.) Suppose $R(\cdot)$ is a $\mu$-strongly convex function with respect to $\|\cdot\|_2$, then $-D(\cdot)$ is $coordinate-wise \left(\frac{\gamma\max_j \|x_j\|_2^2}{n\mu}+1\right)$-smooth relative to $\hh$, and 

(2.) $-D(\cdot)$ is $1-$strongly convex relative to $\hh$.

\end{prop}

{\bf Proof.}
{\em (1.)} Consider $w_1$ and $w_2$ such that $w_2=w_1+te_j$ for some coordinate $j$, namely $w_1$ and $w_2$ only differ in one coordinate. It follows from Proposition \ref{prop:gradient_conjugate} that $R^*(\cdot)$ is $\frac{1}{\mu}$-smooth with respect to $\|\cdot\|_2$, thus we have
\begin{equation*}
\begin{array}{lcl}
R^*\left(-\tfrac{1}{n}X^T w_2\right) & \le & R^*\left(-\tfrac{1}{n}X^T w_1\right) + \left\langle \nabla R^*\left(-\tfrac{1}{n}X^T w_1\right), -\tfrac{1}{n}X^T (w_2 - w_1) \right\rangle + \tfrac{1}{2\mu} \left\|\tfrac{1}{n}X^T(w_2-w_1) \right\|_2^2 \\ \\
& = & R^*\left(-\tfrac{1}{n}X^T w_1\right) + \left\langle -\tfrac{1}{n}X \nabla R^*\left(-\tfrac{1}{n}X^T w_1\right), w_2 - w_1 \right\rangle + \tfrac{t^2}{2n^2\mu} \left\|x_j \right\|_2^2 \\ \\
& \le & R^*\left(-\tfrac{1}{n}X^T w_1\right) + \left\langle -\tfrac{1}{n}X \nabla R^*\left(-\tfrac{1}{n}X^T w_1\right), w_2 - w_1 \right\rangle + \tfrac{\gamma\left\|x_j \right\|_2^2}{n\mu}D_h(w_2, w_1) \ ,
\end{array}
\end{equation*}

where the first inequality follows from smoothness, the equality is from $w_2=w_1+te_j$, and the last inequality utilizes the fact that $h(\cdot)$ is $(\frac{1}{n\gamma})$-strongly convexity with respect to $\|\cdot\|_2$. Therefore it holds that $\hat f(w) := R^*\left(-\tfrac{1}{n}X^T w\right)$ is coordinate-wise $(\frac{\gamma\max_j \|x_j\|_2^2}{n\mu})$-smooth relative to $\hh$.  The proof is completed by noticing that $-D(w)=R^*\left(-\tfrac{1}{n}X^T w\right)+h(w)$.

{\em (2.)} This follows from the additivity property of relative strong convexity (Proposition 1.2 in \cite{lu2018relatively}), whereby $D(\cdot)$ is $1$-strongly concave relative to $h(w)$.
\qed

$ \ $


%



The following theorem states a linear convergence guarantee in the case when $R(\cdot)$ is strongly convex.\medskip

\begin{thm}\label{thm:strong} Suppose $D(\cdot)$ is coordinate-wise $\sigma$-smooth relative to $h(\cdot)$. Consider the Stochastic Generalized Frank-Wolfe method (Algorithm \ref{al:SFW}) or the Randomized Dual Coordinate Mirror Descent method (Algorithm \ref{al:RCMD}), with step-size sequences $\eta_i=\frac{1}{\sigma}$ and $\alpha_i=\frac{n^{-1}\sigma^i}{\sigma^{i+1}-(\sigma-1/n)^{i+1}}.$ 
Under Assumption \ref{ass:smooth} it holds for all $k \ge 1$ that
\begin{equation}\label{eq:thm_strong}
\EE_{\xi_k} \left[P(\bbeta^{k-1})-D(w^{k})\right]\le \frac{D_{\max}}{\left(1+\frac{1}{n\sigma-1}\right)^{k}-1}\le \frac{\gamma M^2}{\left(1+\frac{1}{n\sigma-1}\right)^{k}-1} \ . 
\end{equation}\qed
\end{thm}

Notice that the first inequality in \eqref{eq:thm_strong} shows linear convergence; indeed, in this case it holds that \begin{equation}\label{eq:linear_conv}\frac{ 1}{\left(1+\frac{1}{n\sigma-1} \right)^k-1} \le  n\sigma\left( 1-\tfrac{1}{n\sigma}\right)^k  \ . \end{equation}  (This inequality holds trivially for $k=1$, and induction on $k$ establishes the result for $k \ge 2$.)  Furthermore, when $k$ is large the $-1$  term in the denominator of the left-hand side can be ignored which yields the asymptotic bound $\left( 1-\tfrac{1}{n\sigma}\right)^k D_{\max}$.  The next corollary states the implication of this linear convergence bound in terms of the values $\gamma$ and $\mu$ 
of the $\gamma$-smoothness of $l_1(\cdot), \ldots, l_n(\cdot)$ and the $\mu$-strong convexity of $R(\cdot)$.\medskip

\begin{cor}\label{cor:linear}
Choose $\sigma=\frac{\gamma\max_j \|x_j\|_2^2}{n\mu}+1$ as per Proposition \ref{prop:coordinate smoothness}.  Then Theorem \ref{thm:strong} and \eqref{eq:linear_conv} imply
$$
\EE_{\xi_k} \left[P(\bbeta^{k-1})-D(w^{k})\right]\le \frac{D_{\max}}{\left(1+\frac{1}{\frac{\gamma\max_{j}\|x_j\|_2^2}{\mu} + n-1}\right)^{k}-1}\le D_{\max}\left(\frac{\gamma\max_{j}\|x_j\|_2^2}{\mu}+n\right)\left(1-\frac{1}{n+\frac{\gamma\max_{j}\|x_j\|_2^2}{\mu}}\right)^k.
$$
\end{cor}\medskip

\begin{rem}
{\color{black} Corollary \ref{cor:linear} shows that GSFW requires $O\left((n+\frac{\gamma}{\mu})\log(\frac{1}{\varepsilon})\right)$ iterations to compute an absolute $\varepsilon$-optimal solution of the empirical risk minimization problem with linear prediction \eqref{eq:primal}. This is the same order of convergence rate as SDCA \cite{shalev2013stochastic}.
}
\end{rem}

Before proving Theorem \ref{thm:strong}, we first present an elementary proposition for a separable reference function $h(\cdot)$, whose proof is given in Appendix \ref{app:prop}.\medskip

\begin{prop}\label{prop:separable-expectation}
Suppose $h(\cdot): \ \mathbb{R}^n \rightarrow \mathbb{R}$ is a separable function. Let $j\sim {\cal U}[1,\ldots,n]$. For given $x$, $a$, $y\in \RR^{n}$, define the random
variable $b\in R^{n}$ such that $b_{j}=a_{j}$, and $b_{i}=x_{i}$
for all $i\neq j$. Then:
\[
D_{h}\left(y,a\right)-D_{h}\left(y,x\right)=n\EE_{j}\left(D_{h}\left(y,b\right)-D_{h}\left(y,x\right)\right). \ \qed
\] 
\end{prop}\medskip

We also will use the following proposition whose proof follows easily by induction on $k$.\medskip

\begin{prop}\label{prop:step-series}  Consider the series $\{ \alpha_i \}$ defined by $\alpha_i=\frac{n^{-1}\sigma^i}{\sigma^{i+1}-(\sigma-1/n)^{i+1}}$ for $i \ge 0$, and define the series $\{\bbeta^i\}$ by $\bbeta^{-1} =0$ and $\bbeta^i = (1-\alpha_i)\bbeta^{i-1} + \alpha_i \tbeta^i$ for $i\ge 0$.  Also define $\gamma_i = \left(\frac{n\sigma}{n\sigma-1}\right)^i$ for $i \ge 0$.  Then
$$ \bbeta^k = \frac{\sum_{i=0}^k \gamma_i\tbeta^i}{\sum_{i=0}^k \gamma_i} \ \ \mbox{for~all~} k \ge 0 \ . $$ \qed
\end{prop} \medskip

\textbf{Proof of Theorem \ref{thm:strong}.}

Notice that $\tg^i=\nabla_{j_i}D(w^i)e_{j_i}$, and $w^{i+1}$ is a coordinate update from $w^i$, whereby we have
$$-D(w^{i+1})\le-D(w^{i})-\left\langle  \tg^i,w^{i+1}-w^{i}\right\rangle +\sigma D_{h}(w^{i+1},w^{i})\le-D(w^{i}) \ ,$$ and hence the dual function value sequence $\left\{ D(w^{i})\right\} $ is non-decreasing.

Define $r^{i+1}:=\arg\min_w \left\{ \left\langle -\nabla D(w^{i}),w-w^{i}\right\rangle +\sigma D_{h}(w,w^{i})\right\} $, then we have

\begin{equation}\begin{array}{rl}\label{eq:first_chain_strong}
\EE_{j_i}[-D(w^{i+1})] & \le \EE_{j_i}[ -D(w^{i})-\left\langle \nabla D(w^{i}),w^{i+1}-w^{i}\right\rangle +\sigma D_{h}(w^{i+1},w^{i})]\\ \\
 & = \EE_{j_i}[-D(w^{i})-\tfrac{1}{n}\left(\left\langle \nabla D(w^{i}),r^{i+1}-w^{i}\right\rangle +\sigma D_{h}(r^{i+1},w^{i})\right)]\\ \\
 & \le \EE_{j_i}[-D(w^{i})-\tfrac{1}{n}\left(\left\langle \nabla D(w^{i}),w-w^{i}\right\rangle +\sigma D_{h}(w,w^{i})-\sigma D_{h}(w,r^{i+1})\right)]\\ \\
 & =\EE_{j_i}[-D(w^{i})-\tfrac{1}{n}\left\langle \nabla D(w^{i}),w-w^{i}\right\rangle +\sigma D_{h}(w,w^{i})-\sigma D_{h}(w,w^{i+1})]\\ \\
 & =\EE_{j_i}[-\tfrac{n-1}{n}D(w^{i})-\tfrac{1}{n}\left(D(w^{i})+\left\langle \nabla D(w^{i}),w-w^{i}\right\rangle \right)+\sigma D_{h}(w,w^{i})-\sigma D_{h}(w,w^{i+1})] \ , 
\end{array}\end{equation}
where the first inequality is from the coordinate-wise $\sigma$-smoothness of $D(w)$ relative to $h(w)$ and the fact that $w^{i+1}$ is a coordinate update from $w^i$, the first equality is due to expectation and the separability of $\hh$, the second inequality uses the three-point property (Lemma \ref{tpp}), the second equality uses Proposition \ref{prop:separable-expectation}, and the third equality is just arithmetic rearrangement.

Notice that
\begin{equation}\label{eq:strong_middle}
\begin{array}{cl}
&-D(w^{i})-\left\langle \nabla D(w^{i}),w-w^{i}\right\rangle +n\sigma D_{h}(w,w^{i})-n\sigma D_{h}(w,w^{i+1})\\ \\
 =& -D(w^{i})-\left\langle \nabla_w \phi(\tbeta^i, w^{i}),w-w^{i}\right\rangle +n\sigma D_{h}(w,w^{i})-n\sigma D_{h}(w,w^{i+1})\nonumber \\ \\
 =& -\phi(\tbeta^{i},w)+(n\sigma-1)D_{h}(w,w^{i})-n\sigma D_{h}(w,w^{i+1}) \ ,
\end{array}
\end{equation}

where the last equality utilizes Proposition \ref{prop:equality}. We can then rewrite \eqref{eq:first_chain_strong} (after multiplying by $n$ on both sides) as
\begin{equation}
\EE_{j_i}[-nD(w^{i+1})] \le \EE_{j_i}\left[-(n-1)D(w^{i})-\phi(\tbeta^{i},w)+(n\sigma -1) D_{h}(w,w^{i})-n\sigma D_{h}(w,w^{i+1})\right] \ .
\end{equation}

Multiplying (\ref{eq:strong_middle}) by $\left(\frac{n\sigma}{n\sigma-1}\right)^{i+1}$
and summing over $i=0,\ldots,k-1$, we obtain:
\begin{equation*}
\begin{array}{cl}
& \EE_{\xi_{k}}\left[-\sum_{i=1}^{k}n\left(\frac{n\sigma}{n\sigma-1}\right)^{i}D(w^{i})\right] \\  \\
 \le& \EE_{\xi_{k}}\left[-\sum_{i=1}^{k}(n-1)\left(\frac{n\sigma}{n\sigma-1}\right)^{i}D(w^{i-1})-\sum_{i=1}^{k}\left(\frac{n\sigma}{n\sigma-1}\right)^{i}\phi(\tbeta^{i-1},w)+n\sigma D_{h}(w,w^{0})\right]\\ \\
  \le &\EE_{\xi_{k}}\left[-\sum_{i=1}^{k}(n-1)\left(\frac{n\sigma}{n\sigma-1}\right)^{i}D(w^{i-1})-\left(\sum_{i=1}^{k}\left(\frac{n\sigma}{n\sigma-1}\right)^{i}\right)\phi(\bar{\beta}^{k-1},w)+n\sigma D_{h}(w,w^{0})\right],
\end{array}
\end{equation*}

where the last inequality is from Proposition \ref{prop:step-series} and the convexity of $\phi(\beta,w)$
in $\beta$. Since the sequence $\left\{ D(w^{i})\right\} $ is non-decreasing in $i$ it follows that:
\begin{equation}\label{hello}
\EE_{\xi_k}\left[-\left(\sum_{i=1}^{k}\left(\frac{n\sigma}{n\sigma-1}\right)^{i}\right)D(w^k)\right] \le \EE_{\xi_{k}}\left[-\left(\sum_{i=1}^{k}\left(\frac{n\sigma}{n\sigma-1}\right)^{i}\right)\phi(\bar{\beta}^{k-1},w)\right]+n\sigma D_{h}(w,w^{0}) \ . 
\end{equation}

Let us substitute the following value of $w$ in \eqref{hello}: $w \leftarrow \hw^{k-1}:= \arg\max_w\{\phi(\bbeta^{k-1}, w)\}$, 
which yields:
\[
\left(\sum_{i=1}^{k}\left(\frac{n\sigma}{n\sigma-1}\right)^{i}\right)\EE_{\xi_k}\left[\phi(\bar{\beta}^{k-1},\hw^{k-1})-D({w}^{k})\right]\le n\sigma D_{h}(\hw^{k-1},w^{0})\le n\sigma D_{\max} \ , 
\]
where the last inequality above comes from the definition of $D_{\max}$.
Therefore we have
\[
\EE_{\xi_k}\left[P(\bbeta^{k-1})-D(w^{k})\right]\le\frac{n\sigma}{\left(\sum_{i=1}^{k}\left(\frac{n\sigma}{n\sigma-1}\right)^{i}\right)}D_{\max}=\frac{D_{\max}}{\left(1+\frac{1}{n\sigma-1}\right)^{k}-1},
\]
which furnishes the proof by utilizing Proposition \ref{prop:step-series}.
\qed\medskip

\begin{rem}\label{wilderness}
	Algorithm \ref{al:SFW} as well as the convergence analysis in Theorem \ref{thm:non-strong} and Theorem \ref{thm:strong} can be directly extended to the mini-batch setting. The algorithm extension is accomplished by replacing the single randomly chosen index $j_i$ in the statement of Algorithm \ref{al:SFW} by a random subset of the indices; and the analysis of the algorithm needs to then use the mini-batch of samples instead of a single sample.  More specifically, in the mini-batch version of Algorithm \ref{al:SFW}, we pre-set the batch size $b$, and at iteration $i$ we choose a random subset $B_i\subseteq \{1,2,\ldots,n\}$ of the indices uniformly without replacement such that $|B_i|=b$ at the $i$-th iteration. We then update the predicted values $s$ using only the indices in the subset $B_i$, namely $s_{j}^{i+1}=(1-\eta_{i})s_{j}^{i}+\eta_{i}(x_{j}^{T}\tbeta^{i})$ for $j\in B_i$,  and $s_{j}^{i+1} = s_{j}^{i}$ for $j \not\in B_i$. The update of the substitute gradient becomes $d^{i+1}= \frac{1}{n}X^T\nabla L(s^{i+1}) =  d^i+\frac{1}{n}\sum_{j\in B_i}\left(\dot{l}_{j}(s^{i+1}_{j})-\dot{l}_{j}(s^i_{j})\right)x_{j}$. The rest of the algorithm remains the same as stated in Algorithm \ref{al:SFW}.

By a similar analysis to that of Theorem \ref{thm:non-strong} and Theorem \ref{thm:strong}, we can obtain similar convergence guarantees for the above mini-batch version of GSFW. Essentially we just need to replace $n$ by $n/b$ in the statement of Theorem \ref{thm:non-strong} and Theorem \ref{thm:strong} as well as most of the places in the analysis (except the third and fourth inequalities in \eqref{eq:nonsmooth-theonehand}), after noticing that Corollary \ref{cor:bound} becomes $\|\tg^i\|_2\le \frac{2M\sqrt{b}}{n}$. In particular, for the non-strongly convex case, the convergence guarantee becomes
	$$
	\EE_{\xi_k}\left[ P(\bbeta^k) - D(\bw^k)\right]  \ \le \ \frac{8(\tfrac{n}{b})\gamma M^{2}}{\left(4(\tfrac{n}{b})+k\right)}+\frac{2(\tfrac{n}{b})(2(\tfrac{n}{b})-1)\gamma M^2}{(4(\tfrac{n}{b})+k)(k+1)} \ ,
	$$
	and for the strongly-convex case, the convergence guarantee becomes
	$$
	\EE_{\xi_k} \left[P(\bbeta^{k-1})-D(w^{k})\right]\le \frac{D_{\max}}{\left(1+\frac{1}{\sigma(\tfrac{n}{b})-1}\right)^{k}-1}\le \frac{\gamma M^2}{\left(1+\frac{1}{\sigma(\tfrac{n}{b})-1}\right)^{k}-1} \ . 
	$$
	
Moreover, the updates of Algorithm \ref{al:SFW} in the mini-batch setting can be implemented in parallel as a result of the separability of samples in Algorithm \ref{al:SFW}.
\end{rem}\medskip

\begin{rem}
A natural question to ask next is whether one can achieve an accelerated convergence rate when $R(\cdot)$ is strongly convex, similar to that in \cite{shalev2014accelerated}, \cite{lin2015accelerated}. The answer actually is yes, as one can utilize similar proof techniques as those developed in \cite{lin2015accelerated}. However, the accelerated version may not have a natural interpretation in the primal variables.
\end{rem}

\section{Computational Experiments and Comparisons}
In this section we present the results of some basic numerical experiments where we compare Algorithm \ref{al:SFW} (GSFW) with the following four other stochastic Frank-Wolfe methods in the recent literature:
\begin{itemize}
	\item SCGS -- stochastic gradient sliding algorithm proposed in \cite{lan2016conditional};
	\item SFW -- stochastic Frank-Wolfe algorithm proposed in \cite{hazan2016variance};
	\item SVRF -- stochastic variance reduction Frank-Wolfe algorithm proposed in \cite{hazan2016variance}; and
	\item SCGM -- stochastic conditional gradient method proposed in \cite{mokhtari2018stochastic}.
\end{itemize}

We analyzed the performance of these five algorithms on instances of the following $\ell_1$  norm constrained sparse logistic regression problem:
\begin{equation}
\begin{array}{cl}
\min_{\beta \in \RR^{p}} & P(\beta) = \tfrac{1}{n}\sum_{j=1}^n \ln (1 + \exp(-y_jx_j^T \beta)) \\ \\
\mathrm{s.t.} & \|\beta\|_1 \le \delta \ .
\end{array}
\label{eq:l1_logistic}
\end{equation}
We ran the five stochastic Frank-Wolfe algorithms on ten dataset instances of the constrained logistic regression problem \eqref{eq:l1_logistic} in LIBSVM \cite{chang2011libsvm}.  Here we report on four of these dataset instances, namely a9a, w8a, mushrooms, and gisette, as the results on these four datasets are typical of the results of the other datasets.  Table \ref{tab:stats} describes the dimensions for these four data instances.  We set $\delta=5$ for our experiments with these four datasets, which resulted in solutions on the boundary of the feasible region in all instances.

\begin{table}[]
\centering
\begin{tabular}{|c||c|c|}
\hline
\textbf{dataset}   & \textbf{sample size ($n$)}      & \textbf{feature size ($p$)}  \\ \hline
a9a                & 32561  & 123   \\ \hline
w8a              & 49749  & 300  \\ \hline
mushrooms               & 8124  & 112   \\ \hline
gisette               & 6000  & 5000   \\ \hline
\end{tabular}
\caption{\small{Dimensions and value of $\delta$ for four LIBSVM data instances.}}
\label{tab:stats}
\end{table}

Instead of using one single sample (batch-size equal to $1$) per iteration, we found it far more efficient to run GSFW and SCGM using a mini-batch.  So as not to over-engineer our analysis or unduly bias our results, we used a mini-batch size of $1\%$ of the training data size ($.01n$) and chose all batches randomly without replacement at all iterations (see Remark \ref{wilderness} for a discussion of how to modify GSFW (Algorithm \ref{al:SFW}) using mini-batches, with associated modifications of the computational guarantees).  Note that in theory, algorithms SCGS, SFW, and SVRF all require increasing the batch size to $O(k^3)$, $O(k^2)$, and $O(k)$ at iteration $k$, respectively. In order to retain actual stochasticity of these three methods, we set the maximum batch-size to be $50\%$ of all samples for these three methods; thus when the batch size specified in each of these algorithms is larger than $n/2$, we randomly select $n/2$ samples without replacement in constructing the stochastic gradient estimator.  Table \ref{tab:batch_size} summarizes the above discussion of batch-size modifications for our computational experiments.

\begin{table}[]
\centering
\begin{tabular}{|l||c|c|c|c|c|}\hline
 $ \ \ \ \ \ \ \ \ \ \ $ Batch Size	& GSFW    & SFW  & SCGM    & SCGS  & SVRF \\ \hline\hline
Size consistent with theory & $0.01n$ & $O(k^2)$   & $0.01n$ & $O(k^3)$  & $O(k)$ \\ \hline 
Size used in practice & $0.01n$ & $\min\{k^2,n/2\}$ & $0.01n$ & $\min\{k^3,n/2\}$ & $\min\{k,n/2\}$ \\ \hline
\end{tabular}
\caption{Batch sizes used in the five stochastic Frank-Wolfe methods implemented in our computational experiments.}
\label{tab:batch_size}
\end{table}

\begin{figure}[t!]
\begin{subfigure}[t]{0.45\columnwidth}
\includegraphics[width=\linewidth]{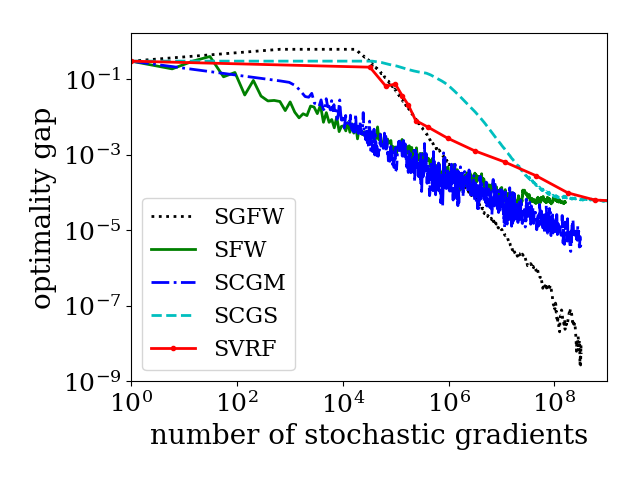}\caption{a9a dataset}
\end{subfigure}
\begin{subfigure}[t]{0.45\columnwidth}
\includegraphics[width=\linewidth]{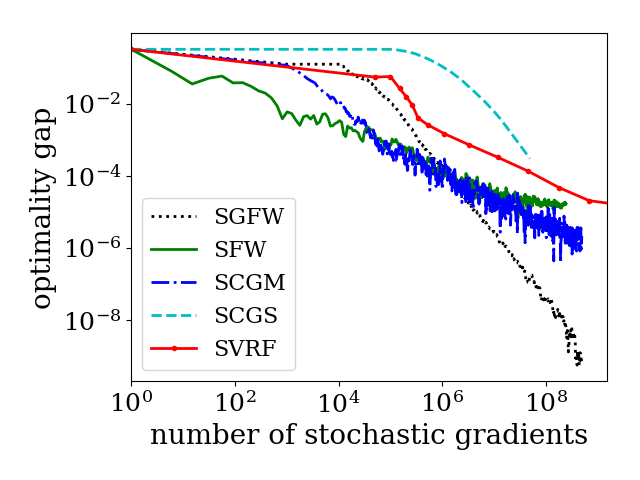}\caption{w8a dataset}  
\end{subfigure}
\\
\begin{subfigure}[t]{0.45\columnwidth}
\includegraphics[width=\linewidth]{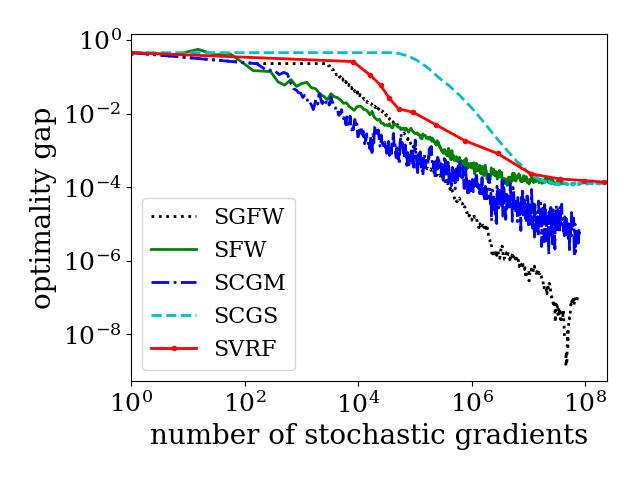}\caption{mushrooms dataset}
\end{subfigure}
\begin{subfigure}[t]{0.45\columnwidth}
\includegraphics[width=\linewidth]{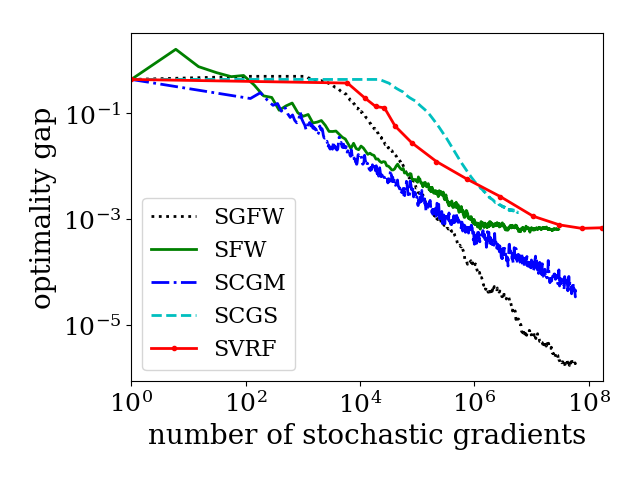}\caption{gisette dataset}
\end{subfigure}
\caption{Figure showing the optimality gap versus the number of stochastic gradient computations (counting one for each sample in each batch) for each of the five different stochastic Frank-Wolfe methods, for the a9a, w8a, mushrooms, and gisette dataset instances.}\label{fig:l1_logistic_sg}
\end{figure}

Figure \ref{fig:l1_logistic_sg} shows the optimality gap versus the number of stochastic gradient computations (counting one for each sample in each batch) for the five different stochastic Frank-Wolfe methods.  In each sub-figure, the vertical axis is the objective value optimality gap $P(\beta^k)-P(\beta^*)$ in $\log$ scale, where $P(\beta^*)$ is estimated after-the-fact using the best solution obtained over the iterations; the horizontal axis in each sub-figure is the number of stochastic gradient computations computed so far, also in $\log$ scale.  Here we see that although the optimality gap of SFW and SCGM may have faster decay initially \textcolor{black}{(perhaps due to better problem-specific constants including their lack of explicit dependence on $n$), nevertheless GSFW decays much faster than the other stochastic Frank-Wolfe variants after a while.  This is consistent with the complexity bound in Theorem \ref{thm:non-strong} that GSFW has a superior order of convergence rate dependence on $\varepsilon$, namely $O(1/\varepsilon)$.}



\begin{figure}[t!]
	\begin{subfigure}[t]{0.45\columnwidth}
		\includegraphics[width=\linewidth]{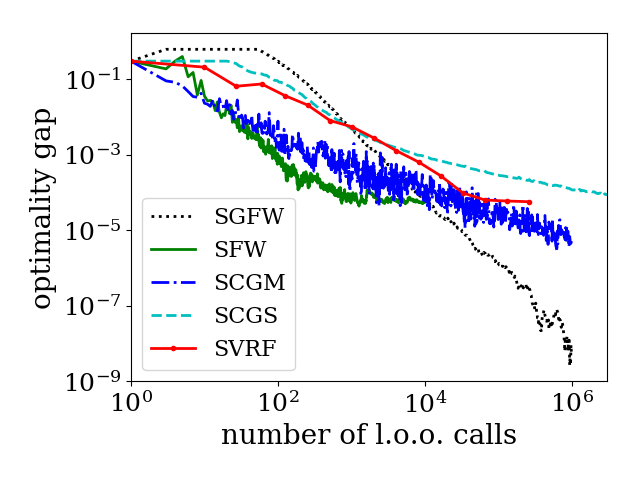}\caption{a9a dataset}
	\end{subfigure}
	\begin{subfigure}[t]{0.45\columnwidth}
		\includegraphics[width=\linewidth]{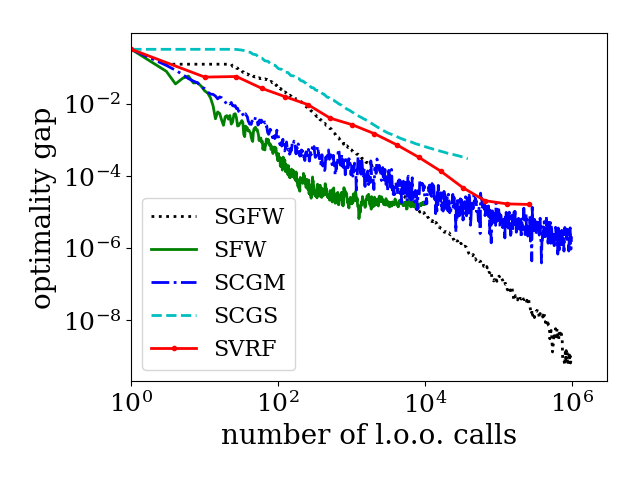}\caption{w8a dataset}  
	\end{subfigure}
	\\
	\begin{subfigure}[t]{0.45\columnwidth}
		\includegraphics[width=\linewidth]{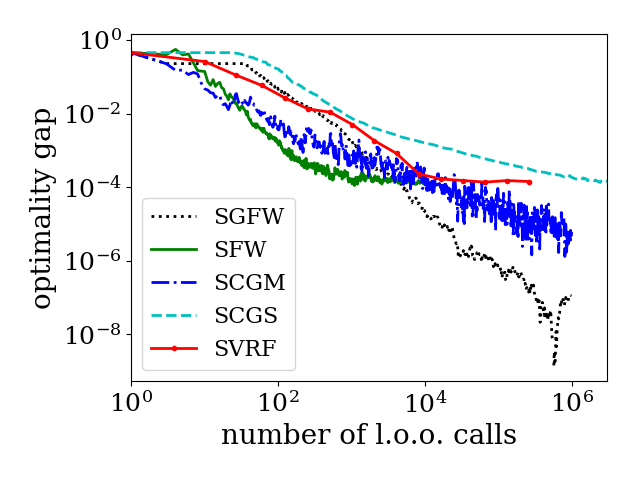}\caption{mushrooms dataset}
	\end{subfigure}
	\begin{subfigure}[t]{0.45\columnwidth}
		\includegraphics[width=\linewidth]{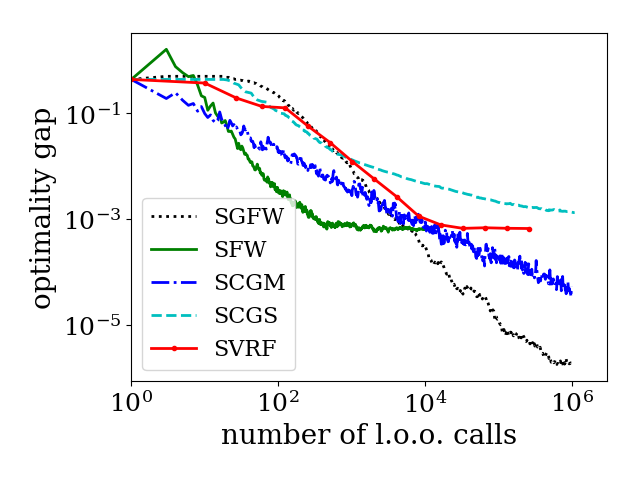}\caption{gisette dataset}
	\end{subfigure}
	\caption{Figure showing the optimality gap versus the number of linear optimization oracle calls for each of the five different stochastic Frank-Wolfe methods, for the a9a, w8a, mushrooms, and gisette dataset instances.}
	\label{fig:l1_logistic_lo}
\end{figure}

Figure \ref{fig:l1_logistic_lo} shows the optimality gap versus the number of linear optimization oracle calls (solving the linear optimization subproblem) for the five different stochastic Frank-Wolfe methods. In each sub-figure, the vertical axis is the objective value optimality gap in $\log$ scale, and the horizontal axis is the number of linear optimization oracle calls, also in $\log$ scale.  Similar to Figure \ref{fig:l1_logistic_sg}, here we see again that the optimality gap of GSFW decays much faster than the other stochastic Frank-Wolfe methods after a while, which is again consistent with the superior order of convergence rate of GSFW dependence on $\varepsilon$, namely $O(1/\varepsilon)$.

Table \ref{tab: comparison} shows the number of stochastic gradient computations and the number of linear optimization oracle calls to achieve an accuracy of $10^{-5}$ for the four LIBSVM datasets evaluated in detail herein, for the (deterministic) Frank-Wolfe method, GSFW, and SCGM.  (SCGS, SFW,  and SVRF do not achieve an accuracy of $10^{-5}$ in a reasonable number of iterations, so we do not present results for these three methods in Table \ref{tab: comparison}.)  Table \ref{tab: comparison} shows that GSFW usually requires fewer gradient computations, while the deterministic Frank-Wolfe method clearly dominates the stochastic Frank-Wolfe methods in term of the number of linear optimization oracle calls. Of course, the Frank-Wolfe method utilizes the exact gradient at each linear optimization oracle call, while the stochastic Frank-Wolfe methods utilize an inexact gradient at each such linear optimization oracle call, and this is likely the reason why the deterministic Frank-Wolfe method requires dramatically fewer linear optimization oracle calls to achieve the given desired optimality gap accuracy. 


\begin{table}[H]
	\centering
	{\small
		\begin{tabular}{|c||c|c|c|c|c|c|}\hline
			& \multicolumn{2}{c|}{Frank-Wolfe} & \multicolumn{2}{c|}{GSFW} & \multicolumn{2}{c|}{SCGM} \\ \hline
			& total sample  & linear  & Stochastic  & linear & Stochastic  & linear  \\ 
			& gradient  & optimization & gradient  & optimization & Gradient  & optimization \\ 
			Dataset & calls ($\times 10^6$) &   oracle calls & calls ($\times 10^6$) & oracle calls & calls ($\times 10^6$) & oracle calls \\ \hline
			a9a      & 14.5  & \bf 448      & \bf 10.3 & 31,900   & 65.3 & 201,000  \\ \hline
			w8a      & 9.65  & \bf 194      & \bf 4.21& 8,470    & 6.92  & 13,900   \\ \hline
			mushroom & 6.44  & \bf 793 & \bf 1.27   & 15,700   & 7.4     & 91,400   \\ \hline
			gisette  & \bf 6.27  & \bf 1045     & 6.56    & 109,000  & --    & --   \\ \hline
		\end{tabular}
		\caption{Comparison of the number of stochastic gradient computations and the number of linear optimization oracle calls, to achieve an accuracy of $10^{-5}$ for solving \eqref{eq:l1_logistic} for the (deterministic) Frank-Wolfe method, GSFW, and SCGM, for the four LIBSVM datasets a9a, w8a, mushrooms, and gisette.} \label{tab: comparison}
	}
	
\end{table}

\appendix

\renewcommand\thesection{\arabic{section}}
\setcounter{section}{-1}
\section*{Appendix}
\renewcommand{\thesection}{A}

\subsection{Properties of Conjugate Functions} \label{app:conjugate}

Recall the definition of the conjugate of a function $\ff$:
$$
f^*(y) := \sup_{x\in \text{dom } f} \{y^T x - f(x)\} \ .
$$
The following properties of conjugate functions are used in this paper:\medskip

\begin{prop}\label{prop:gradient_conjugate}(see \cite{avriel}, \cite{zalinescu}, \cite{kakade2012regularization})
If $f(\cdot)$ is a closed convex function, then ${f}^{**}(\cdot) = f(\cdot)$.  Furthermore:
\begin{enumerate}
\item $\ff$ is $\gamma$-smooth with domain $\RR^p$ with respect to the norm $\|\cdot\|$ if and only if  $f^*(\cdot)$ is $1/\gamma$-strongly convex with respect to the (dual) norm $\|\cdot\|^*$ .
\item If $\ff$ is differentiable and strictly convex, then the following three conditions are equivalent:
\begin{enumerate}
\item $y = \nabla f(x)$
\item $x = \nabla f^*(y)$, and
\item $x^T y = f(x)+f^*(y)$ .
\end{enumerate}
\end{enumerate}
\end{prop}

\subsection{Three-Point Property}

We state here the ``three-point property'' as memorialized by Tseng \cite{tseng}:\medskip\medskip

\begin{lem}\label{tpp} {\bf (Three-Point Property \cite{tseng})} Let $\phi(x)$ be a convex function, and let $\Dh(\cdot, \cdot)$ be the Bregman distance for $\hh$. For a given vector $z$, let
\begin{equation*}
z^+ := \arg\min_{x\in Q} \left\{ \phi(x) + \Dh(x,z) \right\} \ .
\end{equation*}
Then
\begin{equation*}
\phi(x) + \Dh(x,z) \ge \phi(z^+) + \Dh(z^+, z) + \Dh(x,z^+)\ \   for\ all \ x\in Q \ . \qed
\end{equation*}
\end{lem}\medskip\medskip

\subsection{Proof of Proposition \ref{prop:separable-expectation} }\label{app:prop}

Note that
\begin{eqnarray*}
\left\langle \nabla h\left(a\right)-\nabla h\left(x\right),y\right\rangle  & = & \sum_{i=1}^{n}\left\langle \nabla h_{i}\left(a_{i}\right)-\nabla h_{i}\left(x_{i}\right),y_{i}\right\rangle \\
 & = & n\EE_{j}\left\langle \nabla h_{j}\left(a_{j}\right)-\nabla h_{j}\left(x_{j}\right),y_{j}\right\rangle \\
 & = & n\EE_{j}\left\langle \nabla h_{j}\left(b_{j}\right)-\nabla h_{j}\left(x_{j}\right),y_{j}\right\rangle \\
 & = & n\EE_{j}\sum_{i=1}^{n}\left\langle \nabla h_{i}\left(b_{i}\right)-\nabla h_{i}\left(x_{i}\right),y_{i}\right\rangle \\
 & = & n\EE_{j}\left\langle \nabla h\left(b\right)-\nabla h\left(x\right),y\right\rangle\ ,
\end{eqnarray*}
where the second equation is from expectation, and the third and
fourth equation follow because $b_{j}=a_{j}$ and $b_{i}=x_{i}$ for
all $i\neq j$. Using similar logic it also holds that 
\begin{eqnarray*}
\left\langle \nabla h\left(a\right),a\right\rangle -\left\langle \nabla h\left(x\right),x\right\rangle  & = &  n\EE_{j}\left(\left\langle \nabla h\left(b\right),b\right\rangle -\left\langle \nabla h\left(x\right),x\right\rangle\right) \ ,
\end{eqnarray*}
and
\begin{eqnarray*}
h\left(a\right)-h\left(x\right) & = &  n\EE_{j}\left(h\left(b\right)-h\left(x\right)\right).
\end{eqnarray*}
Therefore,
\begin{eqnarray*}
D_{h}\left(y,a\right)-D_{h}\left(y,x\right) & = & \left\langle \nabla h\left(a\right),y-a\right\rangle -\left\langle \nabla h\left(x\right),y-x\right\rangle -\left(h\left(a\right)-h\left(x\right)\right)\\
 & = & \left\langle \nabla h\left(a\right)-\nabla h\left(x\right),y\right\rangle -\left(\left\langle \nabla h\left(a\right),a\right\rangle -\left\langle \nabla h\left(x\right),x\right\rangle \right)-\left(h\left(a\right)-h\left(x\right)\right)\\
 & = & n\EE_{j}\left[\left\langle \nabla h\left(b\right)-\nabla h\left(x\right),y\right\rangle -\left(\left\langle \nabla h\left(b\right),b\right\rangle -\left\langle \nabla h\left(x\right),x\right\rangle \right)-\left(h\left(b\right)-h\left(x\right)\right)\right]\\
 & = & n\EE_{j}\left[\left\langle \nabla h\left(b\right),y-b\right\rangle -\left\langle \nabla h\left(x\right),y-x\right\rangle -\left(h\left(b\right)-h\left(x\right)\right)\right]\\
 & = & n\EE_{j}\left[D_{h}\left(y,b\right)-D_{h}\left(y,x\right)\right] \ .  \  \qed
\end{eqnarray*}\medskip

\subsection{Connections and Comparisons between GSFW and Stochastic Dual Coordinate Ascent Methods}\label{sec:connections-to-shalev}
In this subsection we discuss connections and comparisons between GSFW (Algorithm \ref{al:SFW}, and equivalently Algorithm \ref{al:RCMD}) and SDCA.  The traditional analysis in SDCA \cite{shalev2013stochastic}\cite{lin2015accelerated}\cite{shalev2014accelerated}\cite{shalev2012proximal} is premised on the assumption that $R(\cdot)$ is a strongly convex function, whereby the first term in the dual objective \eqref{eq:dual} is a smooth function.  In contrast, for our GSFW method $R(\cdot)$ need not have any such structure; indeed in the Frank-Wolfe setting $R(\cdot)$ can be an indicator function of the (primal) feasible region whereby the first term in the dual problem \eqref{eq:dual} is then non-differentiable.  Viewing GSFW through the dual (Algorithm \ref{al:RCMD}), we compute a subgradient of the first part of the dual objective by calling a linear optimization oracle in the primal space, and this subgradient is used in the dual mirror descent algorithm for solving \eqref{eq:dual}.  From the SDCA perspective, this is the first such version of SDCA that does not require a strongly convex regularizer $R(\cdot)$.  

(As a thought exercise, it is surely possible to start with a non-strongly convex $R(\cdot)$ and then add a tiny strongly convex regularizer term based on a target optimality tolerance and other continuity parameters [though such tuning can be tricky to do in practice], and then use SDCA rather than GSFW.  However, in this approach the subproblem that needs to be solved at each iteration requires a projection step onto the feasible region as opposed to solving a linear optimization subproblem, which goes beyond and can be much more computationally demanding than Frank-Wolfe in certain settings.  For example, in matrix completion (Example \ref{robespierre}), solving the linear optimization oracle with the nuclear norm ball requires the computation of the largest eigenvector/eigenvalue pair, while doing the projection requires a full eigendecomposition, and therefore can be significantly more computationally burdensome.  Of course, this is one of the key reasons why Frank-Wolfe methods have been so extensively studied in the past decade.)

GSFW has closer connections to SDCA in the case when $R(\cdot)$ is strongly convex -- whereby the first term in the dual objective is smooth.  In this case Randomized Coordinate Mirror Descent (Algorithm \ref{al:RCMD}) can be viewed as a variant of SDCA with a specific updating rule based on the mirror descent methodology. This perspective provides a new interpretation for SDCA in the primal space as a variant of a Frank-Wolfe based method.  In the previous literature involving SDCA, even though one can rewrite SDCA entirely in the primal space \cite{shalev2016sdca}, there are still explicit dual variables which lack intuition or interpretation in the primal space, in contrast to our Algorithm \ref{al:SFW} and Algorithm \ref{al:RCMD} equivalency.  

In Algorithm \ref{al:RCMD} the coordinate update at each iteration requires the solution of a univariate problem of the following form for a suitably given scalar $c_{j_i}$:
\begin{equation}\label{eq:our_update}
\min_{w_{j_i}} c_{j_i}w_{j_i} +l_{j_i}^*(w_{j_i}) \ ,
\end{equation}
in comparison with the basic version SDCA algorithms in \cite{shalev2013stochastic} or \cite{lin2015accelerated} for which the coordinate update at each iteration requires the solution of the following slightly different univariate problem for suitably given scalars $c_{j_i}$ and $b_{j_i}$:
\begin{equation}\label{eq:their_update}
\min_{w_{j_i}} c_{j_i}w_{j_i} +l_{j_i}^*(w_{j_i}) + b_{j_i}w_{j_i}^2\ .
\end{equation}
On the other hand, we point out that there are variants of SDCA that do not require solving \eqref{eq:their_update}, see the update rules [III], [IV] and [V] in \cite{shalev2012proximal}, but these rules are still different from \eqref{eq:our_update}.  Indeed, the unaccelerated version of \cite{lin2015accelerated} can be viewed as a randomized coordinate method with a composite function, and has the following update:
$$
w^{i+1}_{j_i}=\arg\min_{w_{j_i}}\left\{ c_{j_i} w_{j_i} + \tfrac{1}{2\eta}(w_{j_i}-w_{j_i}^i)^{2} + \tfrac{1}{n}l^*_{j_i}(w_{j_i}) \right\} \ ,
$$
where $c_{j_i}$ is one coordinate of the gradient.  The above update can be viewed as a coordinate mirror descent method update with reference function $h(w):=\tfrac{1}{2\eta}\|w\|^2 + \tfrac{1}{n}L^*(w)$. This is similar to the equivalence of composite optimization and mirror descent in the deterministic case discussed in Section 3.3 of \cite{lu2018relatively}.

\subsection{Regarding Randomized Coordinate Mirror Descent with Non-smooth Functions}\label{sec:nonsmooth-cd}
Since the seminal work of Nesterov \cite{nesterov2012efficiency}, there have been many research results on randomized coordinate descent for convex minimization of a general smooth objective function; however, there has not been much research on randomized coordinate descent in the general non-smooth setting.  \cite{nesterov2014subgradient} develops a randomized block-coordinate method for some specially structured problems.  Also, many papers consider a composite objective function $f(\cdot) := \hat f(\cdot) + \tau(\cdot)$ where $\hat f(\cdot)$ is smooth and $\tau(\cdot)$ is non-smooth, separable, and computationally friendly, see in particular \cite{nesterov2012efficiency}, \cite{richtarik2014iteration}, \cite{lu2015complexity} among many others.  These methods cannot be applied to the general non-smooth case so far as we can tell.  We can consider the case where $f(\cdot)$ is the sum of a non-smooth Lipschitz continuous function $\bar f(\cdot)$ and a separable strongly convex function, namely $f(w) := \bar f(w) + \sum_{j=1}^n \tau_j(w_j)$ where $\tau_j(\cdot): \mathbb{R} \rightarrow \mathbb{R}$ is strongly convex on its domain.  Then $f(\cdot)$ is in the format of the dual problem \eqref{eq:dual}, whereby Theorem {\ref{thm:non-strong}} can be applied to show that the randomized coordinate mirror descent method (Algorithm \ref{al:RCMD}) yields the indicated computational guarantees in this case.  This shows that non-smooth convex optimization can be tackled by randomized coordinate descent if the objective function is the sum of a non-smooth function and a strongly convex separable function.



\bibliographystyle{amsplain}
\bibliography{LF-papers}
\end{document}